\documentclass[11pt]{article}

\usepackage{amsmath}
\usepackage{amsthm}
\usepackage{amssymb}
\usepackage{amscd}

\theoremstyle{definition}
\newtheorem{theorem}{Theorem}[section]
\newtheorem{prop}[theorem]{Proposition}

\newtheorem{corollary}[theorem]{Corollary}
\newtheorem{definition}[theorem]{Definition}

\newtheorem{remark}[theorem]{Remark}

\numberwithin{equation}{section}
\newenvironment{demo}[1]{%
  \trivlist
  \item[\hskip\labelsep
        {\it #1.}]
}{%
\hfill\qedsymbol
  \endtrivlist
}

\newcommand\Int{\mathbb{Z}}

\newcommand\Comp{\mathbb{C}}

\newcommand\RSTab{\operatorname{RSTab}} 
\newcommand\CSTab{\operatorname{CSTab}} 

\newcommand\Orth{\mathbf{O}}
\newcommand\SOrth{\mathbf{SO}}
\newcommand\GL{\mathbf{GL}}
\newcommand\Symp{\mathbf{Sp}}
\newcommand\Res{\operatorname{Res}}
\newcommand\LR{\operatorname{LR}}
\newcommand\Par{\mathcal{P}}
\newcommand\Sym{\mathfrak{S}}
\newcommand\sgn{\operatorname{sgn}}
\newcommand\twedge{{\textstyle\bigwedge}} 

\newcommand\trans{{}^t\!}
\renewcommand\tilde{\widetilde}
\newcommand\ep{\varepsilon}

\newcommand\vectx{\boldsymbol{x}}

\allowdisplaybreaks

\title{
Pieri Rules for Classical Groups 
and Equinumeration between Generalized Oscillating Tableaux and Semistandard Tableaux
}

\author{
Soichi OKADA
\footnote{
Graduate School of Mathematics, Nagoya University, 
Furo-cho, Chikusa-ku, Nagoya 464-8602, Japan. 
e-mail: {\tt okada@math.nagoya-u.ac.jp}
}
\footnote{
This work was partially supported by the JSPS Grants-in-Aid for Scientific Research No. 24340003.
}
}

\date{
}

\begin{document}

\maketitle

\begin{abstract}
We present several equinumerous results between generalized oscillating tableaux 
and semistandard tableaux and give a representation-theoretical proof to them.
As one of the key ingredients of the proof, we provide Pieri rules for the symplectic and 
orthogonal groups.
\end{abstract}

\section{%
Introduction
}

Various types of tableaux are studied in combinatorics and representation theory.
Representation theory provides powerful tools to prove theorems concerning tableaux
and helpful clues to generalize them.
The aim of this article is to give a representation-theoretical proof of 
generalizations and variants of refined Burrill conjecture for oscillating tableaux.
These generalizations and variants are equinumerous results between generalized oscillating tableaux 
and semistandard tableaux.

A standard tableau of shape $\lambda$ can be viewed as an increasing sequence 
$\emptyset = \lambda^{(0)} \subset \lambda^{(1)} \subset \dots \subset \lambda^{(k-1)} 
\subset \lambda^{(k)} = \lambda$ of partitions such that 
the diagram of $\lambda^{(i)}$ is obtained from that of $\lambda^{(i-1)}$ 
by adding one cell for each $i$.
The notion of oscillating tableaux is a generalization of standard tableaux.
For a nonnegative integer $k$ and a partition $\lambda$, 
an \emph{oscillating tableau} of length $k$ and shape $\lambda$ is 
a sequence $(\lambda^{(i)})_{i=0}^k = (\lambda^{(0)}, \lambda^{(1)}, \dots, \lambda^{(k)})$ 
of partitions satisfying the following two conditions:
\begin{enumerate}
\item[(i)]
$\lambda^{(0)} = \emptyset$ and $\lambda^{(k)} = \lambda$.
\item[(ii)]
The diagram of $\lambda^{(i)}$ is obtained from that of $\lambda^{(i-1)}$ 
by adding or removing one cell for each $i$.
\end{enumerate}
Krattenthaler \cite{K} and Burrill--Courtiel--Fusy--Melczer--Mishna \cite{BCFMM} 
gave bijective proofs of the following theorem, which is a refinement of 
Burrill's conjecture \cite[Conjecture~6.2.1]{B}.

\begin{theorem}
\label{thm:Burrill}
(\cite[Theorem~3]{K}, \cite[Theorem~1]{BCFMM})
For nonnegative integers $k$, $n$ and $m$, the following two sets are equinumerous:
\begin{enumerate}
\item[(a)]
The set of oscillating tableaux $(\lambda^{(i)})_{i=0}^k$ of length $k$ and shape $(m)$ 
such that $l(\lambda^{(i)}) \le n$ for each $i$.
\item[(b)]
The set of standard tableaux whose shape $\lambda$ satisfies $|\lambda| = k$, 
$l(\lambda) \le 2n$ and $c(\lambda) = m$.
\end{enumerate}
Here $l(\lambda)$ is the length of a partition $\lambda$, 
and $c(\lambda)$ is the number of columns of odd length in the diagram of $\lambda$.
\end{theorem}

In this paper, we present generalizations and variants of Theorem~\ref{thm:Burrill} 
from the view point of representation theory of classical groups.
Our main results in the symplectic group case can be stated as follows.
(See Theorem~\ref{thm:main2} for similar results in the orthogonal group case.)
We remark that Krattenthaler \cite[Theorem~4]{K} gave a bijective proof to part (1) of 
the following theorem.

\begin{theorem}
(Theorem~\ref{thm:main2} (1) and (2) below)
\label{thm:main}
Let $\alpha = (\alpha_1, \dots, \alpha_k)$ be a sequence of nonnegative integers of length $k$, 
and $m$ and $n$ nonnegative integers.
Then we have
\begin{enumerate}
\item[(1)]
The following two sets are equinumerous:
\begin{enumerate}
\item[(a)]
The set of sequences
$$
\emptyset = \lambda^{(0)} \supset \lambda^{(1)} \subset \lambda^{(2)} 
\supset \lambda^{(3)} \subset \lambda^{(4)} \supset \dots 
\subset \lambda^{(2k-2)} \supset \lambda^{(2k-1)} \subset \lambda^{(2k)} = (m)
$$
of partitions with $l(\lambda^{(i)}) \le n$ such that 
$\lambda^{(2i-2)}/\lambda^{(2i-1)}$ and $\lambda^{(2i)}/\lambda^{(2i-1)}$ are horizontal strips 
and that $|\lambda^{(2i-2)}/\lambda^{(2i-1)}| + |\lambda^{(2i)}/\lambda^{(2i-1)}| = \alpha_i$.
\item[(b)]
The set of column-strict tableaux (a.k.a. semistandard tableaux) 
of weight $\alpha$ whose shape $\lambda$ satisfies 
$l(\lambda) \le 2n$ and $c(\lambda) = m$.
\end{enumerate}
\item[(2)]
The following two sets are equinumerous:
\begin{enumerate}
\item[(a)]
The set of sequences
$$
\emptyset = \lambda^{(0)} \subset \lambda^{(1)} \supset \lambda^{(2)} 
\subset \lambda^{(3)} \supset \lambda^{(4)} \subset \dots 
\supset \lambda^{(2k-2)} \subset \lambda^{(2k-1)} \supset \lambda^{(2k)} = (m)
$$
of partitions with $l(\lambda^{(i)}) \le n$ such that 
$\lambda^{(2i-1)}/\lambda^{(2i-2)}$ and $\lambda^{(2i-1)}/\lambda^{(2i)}$ are vertical strips 
and that $|\lambda^{(2i-1)}/\lambda^{(2i-2)}| + |\lambda^{(2i-1)}/\lambda^{(2i)}| = \alpha_i$.
\item[(b)]
The set of row-strict tableaux of weight $\alpha$ whose shape $\lambda$ satisfies 
$l(\lambda) \le 2n$ and $c(\lambda) = m$.
\end{enumerate}
\end{enumerate}
See Sections 2 and 5 for undefined terminologies.
\end{theorem}

Here we give an outline of our proof of Theorem~\ref{thm:main}.
Let $V = \Comp^{2n}$ be the defining representation of the general linear group 
$\GL_{2n} = \GL_{2n}(\Comp)$, which is also the defining representation 
of the symplectic group $\Symp_{2n} = \Symp_{2n}(\Comp)$.
We consider the tensor products
$$
S^\alpha (V) = S^{\alpha_1} (V) \otimes \dots \otimes S^{\alpha_k} (V)
\quad
\text{and}
\quad
\twedge^\alpha (V) = \twedge^{\alpha_1} (V) \otimes \dots \otimes \twedge^{\alpha_k} (V)
$$
of the symmetric and exterior powers of $V$.
We compute in two ways the multiplicities of the irreducible $\Symp_{2n}$-module $V_{\langle (m) \rangle}$ 
corresponding to the one-row partition $(m)$ in $S^\alpha(V)$ and $\twedge^\alpha(V)$.

One of the key ingredients of the proof is the Pieri rules for classical groups, 
which are another contribution of this paper.
The Pieri rule (resp. dual Pieri rule) for $\Symp_{2n}$ describes 
the irreducible decomposition of the tensor product 
of an irreducible representation with $S^r(V)$ (resp. $\twedge^r(V))$.
By iteratively applying the Pieri rule given in 
Theorem~\ref{thm:Pieri} (1) (resp. the dual Pieri rule in Theorem~\ref{thm:dualPieri} (1)), 
we see that 
the multiplicity of $V_{\langle (m) \rangle}$ in $S^\alpha(V)$ (resp. $\twedge^\alpha(V)$) 
is equal to the number of the sequences of partitions described in (a) of Theorem~\ref{thm:main} (1) 
(resp. (2)).

On the other hand, 
by using the classical Pieri rule for $\GL_{2n}$, 
we obtain the $\GL_{2n}$-module decompositions 
\begin{equation}
\label{eq:outline1}
S^\alpha (V) 
\cong 
\bigoplus_{l(\lambda) \le 2n} V_\lambda^{\oplus \CSTab(\lambda, \alpha)},
\quad
\twedge^\alpha (V) 
\cong 
\bigoplus_{l(\lambda) \le 2n} V_\lambda^{\oplus \RSTab(\lambda, \alpha)},
\end{equation}
where $\lambda$ runs over all partitions of length $\le 2n$, 
$V_\lambda$ is the irreducible representation of $\GL_{2n}$ corresponding to $\lambda$, 
and $\CSTab(\lambda,\alpha)$ (resp. $\RSTab(\lambda,\alpha)$) denotes the set of 
column-strict (resp. row-strict) tableaux of shape $\lambda$ and weight $\alpha$.
Another key ingredient of the proof is the following restriction multiplicity formula 
(see Theorem~\ref{thm:rest}):
\begin{equation}
\label{eq:outline2}
\big[
 \Res^{\GL_{2n}}_{\Symp_{2n}} V_\lambda: V_{\langle (m) \rangle}
\big]
 =
\begin{cases}
 1 &\text{if $c(\lambda) = m$,} \\
 0 &\text{otherwise.}
\end{cases}
\end{equation}
It follows from (\ref{eq:outline1}) and (\ref{eq:outline2}) 
that the multiplicity of $V_{\langle (m) \rangle}$ in 
$S^\alpha(V)$ (resp. $\twedge^\alpha(V)$) is equal to the number of tableaux described in 
(b) of Theorem~\ref{thm:main} (1) (resp. (2)).
In this way, we prove Theorem~\ref{thm:main}.

The remaining of this paper is organized as follows.
In Section~2 we review the representation theory of classical groups.
In Sections~3 and 4, we prove the Pieri rules and the dual Pieri rules for the symplectic 
and orthogonal groups.
In Section~5, we present variants/generalizations of Theorem~\ref{thm:Burrill} 
and give a representation-theoretical proof to them.

\section{%
Preliminaries
}

In this section, first we recall some definitions on partitions.
Then we review the representation theory of classical groups 
and collect several facts which will be used in the remaining of the paper.
See \cite{FH}, \cite{GW} and \cite{O} for the representation theory of classical groups.

\subsection{%
Combinatorics of partitions
}

A partition is a weakly decreasing sequence 
$\lambda = (\lambda_1, \lambda_2, \dots)$ of nonnegative integers 
such that $\sum_{i \ge 1} \lambda_i$ is finite. 
The length of $\lambda$, denoted by $l(\lambda)$, is the number of nonzero entries of $\lambda$, 
and the size of $\lambda$, denoted by $|\lambda|$, is the sum of entries of $\lambda$.
A partition $\lambda$ is often identified with its diagram, 
which is a left-justified array of $|\lambda|$ cells with $\lambda_i$ cells in the $i$th row.
The conjugate partition $\lambda'$ of a partition $\lambda$ is the partition 
whose diagram is obtained by reflecting the diagram of $\lambda$ along the main diagonal.
We denote by $\Par$ the set of all partitions.

For two partitions $\lambda$ and $\mu$, we write $\mu \subset \lambda$ 
if $\mu_i \le \lambda_i$ for all $i$.
Then the skew diagram $\lambda/\mu$ is defined to be the set-theoretical difference 
of the diagrams of $\lambda$ and $\mu$.
The size of the skew diagram is defined by $|\lambda/\mu| = |\lambda| - |\mu|$.
We say that the skew diagram $\lambda/\mu$ is a \emph{horizontal $r$-strip} 
if it contains at most one cell in each column and $|\lambda/\mu| = r$.
Dually, we say that the skew diagram $\lambda/\mu$ is a \emph{vertical $r$-strip} 
if it contains at most one cell in each row and $|\lambda/\mu| = r$.
Note that $\lambda/\mu$ is a horizontal strip 
if and only if $\lambda_1 \ge \mu_1 \ge \lambda_2 \ge \mu_2 \ge \dots$.

\subsection{%
Representation theory of $\GL_N$
}

It is well-known that the irreducible polynomial representations of the general linear group $\GL_N
 = \GL_N(\Comp)$ are parametrized by partitions of length $\le N$.
We denote by $V_\lambda$ and $S_\lambda$ 
the irreducible representation and its character corresponding to a partition $\lambda$ respectively.
If $\lambda$ is a one-row partition $(r)$ (resp. a one-column partition $(1^r)$), 
then we write $H_r = S_{(r)}$ (resp. $E_r = S_{(1^r)}$).
Note that $H_r$ (resp. $E_r$) is the characters of the symmetric power $S^r(V)$ 
(resp. the exterior power $\twedge^r(V)$) of the defining representation $V = \Comp^N$ of $\GL_N$.
We also use the notations $V_{\lambda,\GL_N}$, $S_{\lambda,\GL_N}$, $H_{r,\GL_N}$ and $E_{r,\GL_N}$ 
to avoid confusions.
Let $R(\GL_N)$ be the representation ring of $\GL_N$.
Then $R(\GL_N)$ is a free $\Int$-module with basis 
$\{ S_\lambda : \lambda \in \Par, \ l(\lambda) \le N \}$.

Let $\Lambda$ be the ring of symmetric functions. 
Let $s_\lambda \in \Lambda$ be the Schur function associated with a partition $\lambda$ 
and $h_r = s_{(r)} \in \Lambda$ the complete symmetric function of degree $r$.
If $\pi = \pi_{\GL_N} : \Lambda \to R(\GL_N)$ is the ring homomorphism defined by 
$\pi(h_r) = H_r$ for $r \ge 1$, then we have
$$
\pi_{\GL_N}(s_\lambda) = \begin{cases}
 S_\lambda &\text{if $l(\lambda) \le N$,} \\
 0 &\text{otherwise.}
\end{cases}
$$
For partitions $\lambda$, $\mu$ and $\nu$, we denote by $\LR^\lambda_{\mu,\nu}$ 
the Littlewood-Richardson coefficient, which is defined by the relation
$$
s_\mu \cdot s_\nu = \sum_{\lambda \in \Par} \LR^\lambda_{\mu,\nu} s_\lambda.
$$
It is known that $\LR^\lambda_{\mu,\nu} = 0$ unless $|\lambda| = |\mu| + |\nu|$, 
$\lambda \supset \mu$ and $\lambda \supset \nu$.

The following is the classical Pieri rule for $\GL_N$.

\begin{prop}
\label{prop:GL-Pieri}
Let $V$ be the defining representation of $\GL_N$.
For a partition $\mu$ of length $\le N$ and a nonnegative integer $r$, we have
\begin{equation}
\label{eq:GL-Pieri}
V_\mu \otimes S^r(V) \cong \bigoplus_\lambda V_\lambda,
\quad
V_\mu \otimes \twedge^r(V) \cong \bigoplus_\rho V_\rho,
\end{equation}
where $\lambda$ (resp. $\rho$) runs over all partitions of length $\le N$ 
such that $\lambda/\mu$ is a horizontal $r$-strip 
(resp. $\rho/\mu$ is a vertical $r$-strip).
In other words, for partitions $\lambda$ and $\mu$ and a nonnegative integer $r$, we have
\begin{align}
\label{eq:GL-Pieri1}
\LR^\lambda_{\mu,(r)} = \begin{cases}
 1 &\text{if $\lambda/\mu$ is a horizontal $r$-strip,} \\
 0 &\text{otherwise,}
\end{cases}
\\
\label{eq:GL-Pieri2}
\LR^\lambda_{\mu,(1^r)} = \begin{cases}
 1 &\text{if $\lambda/\mu$ is a vertical $r$-strip,} \\
 0 &\text{otherwise.}
\end{cases}
\end{align}
\end{prop}

\subsection{%
Representation theory of $\Symp_{2n}$
}

Next we consider the symplectic group $\Symp_{2n} = \Symp_{2n}(\Comp)$.
The finite-dimensional irreducible representations of $\Symp_{2n}$ 
are indexed by partitions of length $\le n$.
Let $\Par(\Symp_{2n})$ be the set of all partitions of length $\le n$.
We denote by $V_{\langle \lambda \rangle} = V_{\langle \lambda \rangle, \Symp_{2n}}$ 
and $S_{\langle \lambda \rangle} = S_{\langle \lambda \rangle, \Symp_{2n}}$ 
the irreducible representation and its character of $\Symp_{2n}$ corresponding 
to a partition $\lambda$ with $l(\lambda) \le n$.
Let $V = \Comp^{2n}$ be the defining representation of $\Symp_{2n}$, and
denote by $H_r$ (resp. $E_r$) 
the character of $\Symp_{2n}$ on $S^r(V)$ (resp. $\twedge^r(V)$).
Then $S^r(V)$ is the irreducible representation corresponding to the one-row partition $(r)$, 
while $\twedge^r(V)$ is not irreducible if $r \ge 2$ and the quotient $\twedge^r(V)/\twedge^{r-2}(V)$ 
is the irreducible representation corresponding to the one-column partition $(1^r)$.
Let $R(\Symp_{2n})$ be the representation ring of $\Symp_{2n}$ 

For an arbitrary partition $\lambda$, we define 
the corresponding \emph{symplectic Schur function} $s_{\langle \lambda \rangle} \in \Lambda$ by putting
\begin{equation}
\label{eq:Sp-Schur}
s_{\langle \lambda \rangle}
 =
\frac{1}{2}
\det \left(
 h_{\lambda_i-i+j} + h_{\lambda_i-i-j+2}
\right)_{1 \le i, j \le l(\lambda)}.
\end{equation}
Let $\pi_{\Symp_{2n}} : \Lambda \to R(\Symp_{2n})$ be the ring homomorphism defined by 
$\pi_{\Symp_{2n}}(h_r) = H_r$ for $r \ge 1$.
The image of a symplectic Schur function under $\pi_{\Symp_{2n}}$ can be 
expressed as a linear combination of irreducible characters 
(in fact, it is $0$ or an irreducible character up to sign), 
by using the following algorithm (see \cite{King} and \cite{KT}).

\begin{prop}
\label{prop:Sp-algorithm}
Let $\lambda$ be a partition.
\begin{enumerate}
\item[(1)]
If $\lambda \in \Par(\Symp_{2n})$, then we have
$$
\pi_{\Symp_{2n}} ( s_{\langle \lambda \rangle} )
 = S_{\langle \lambda \rangle}.
$$
\item[(2)]
In general, the image $\pi_{\Symp_{2n}} ( s_{\langle \lambda \rangle} )$ is computed as follows.
We put $r = \lambda_1$ and
$$
\alpha = (\lambda'_1, \lambda'_2-1, \dots, \lambda'_r-(r-1)).
$$
\begin{enumerate}
\item[(a)]
If $\alpha$ has an entry larger than or equal to $2n+r+2$, 
then $\pi_{\Symp_{2n}} ( s_{\langle \lambda \rangle} ) = 0$.
\item[(b)]
If $\alpha_i + \alpha_j = 2n+2$ for some $i$ and $j$, 
then $\pi_{\Symp_{2n}} ( s_{\langle \lambda \rangle} ) = 0$.
\item[(c)]
Otherwise, suppose that $\alpha_1 > \dots > \alpha_p \ge n+2 > \alpha_{p+1}$ 
and define a sequence $\beta$ by putting
$$
\beta
 =
(2n+2-\alpha_1, \dots, 2n+2-\alpha_p, \alpha_{p+1}, \dots, \alpha_r).
$$
Let $\gamma$ be the rearrangement of $\beta$ in decreasing order 
and $\sigma \in \Sym_r$ a permutation satisfying $\sigma(\beta) = \gamma$.
If $\mu$ is the partition given by $\gamma = (\mu'_1, \mu'_2-1, \dots, \mu'_r-(r-1))$, 
then $\mu \in \Par(\Symp_{2n})$ and 
$$
\pi_{\Symp_{2n}} ( s_{\langle \lambda \rangle} )
 = (-1)^p \sgn(\sigma) S_{\langle \mu \rangle}.
$$
\end{enumerate}
\item[(3)]
In particular, if $l(\lambda) = n+1$, 
then $\pi_{\Symp_{2n}} ( s_{\langle \lambda \rangle} ) = 0$.
\end{enumerate}
\end{prop}

We have the following relations in the ring $\Lambda$ of symmetric functions, 
from which we can derive identities involving irreducible characters of $\Symp_{2n}$ 
by applying $\pi_{\Symp_{2n}}$.

\begin{prop}
\label{prop:Sp-schur-rel}
\begin{enumerate}
\item[(1)]
(Newell, Littlewood)
For any partitions $\mu$ and $\nu$, we have
\begin{equation}
\label{eq:Sp-LR}
s_{\langle \mu \rangle} \cdot s_{\langle \nu \rangle}
 =
\sum_{\lambda \in \Par}
 \left(
  \sum_{\tau, \, \xi, \, \eta \in \Par} \LR^\mu_{\tau,\xi} \LR^\nu_{\tau,\eta} \LR^\lambda_{\xi,\eta}
 \right)
 s_{\langle \lambda \rangle}.
\end{equation}
\item[(2)]
(Littlewood)
For any partition $\lambda$, we have
\begin{equation}
\label{eq:Sp-GL}
s_\lambda
 =
\sum_{\mu \in \Par}
 \left(
  \sum_{\kappa \in \mathcal{E}'} \LR^\lambda_{\kappa,\mu}
 \right)
 s_{\langle \mu \rangle},
\end{equation}
where $\mathcal{E}'$ is the set of all partitions whose column lengths are all even.
\end{enumerate}
\end{prop}

Let $\{ f^C_r(x) \}$ be the sequence of Laurent polynomials defined by
$$
f^C_r(x)
 =
\begin{cases}
 x^{r+1} - x^{-r-1} &\text{if $r \ge 0$,} \\
 0 &\text{if $r < 0$.}
\end{cases}
$$
For a sequence $\alpha = (\alpha_1, \dots, \alpha_n)$ of integers 
and a sequence $\vectx = (x_1, \dots, x_n)$ of indeterminates, 
we put
$$
A^C_\alpha(\vectx)
 =
\left( f^C_{\alpha_j+n-j}(x_i) \right)_{1 \le i, j \le n}.
$$
Then the ratio $\det A^C_{\alpha}(\vectx)/\det A^C_{\emptyset}(\vectx)$, 
where $\emptyset = (0, \dots, 0)$, is a Laurent polynomial in $x_1, \dots, x_n$.
Then the Weyl character formula is rephrased as follows:

\begin{prop}
\label{prop:Sp-Weyl}
If $X \in \Symp_{2n}$ has the eigenvalues $x_1, \dots, x_n, x_1^{-1}, \dots, x_n^{-1}$, 
then we have
\begin{equation}
\label{eq:Sp-Weyl}
S_{\langle \lambda \rangle}(X)
 =
\frac{ \det A^C_\lambda (\vectx) }
     { \det A^C_\emptyset (\vectx) }.
\end{equation}
\end{prop}

\subsection{%
Representation theory of $\Orth_N$
}

The finite-dimensional irreducible representations of the orthogonal group $\Orth_N = \Orth_N(\Comp)$ 
are parametrized by partitions such that the sum of the lengths of the first two columns 
is at most $N$.
Let $\Par(\Orth_N)$ be the set of partitions $\lambda$ satisfying $\lambda'_1 + \lambda'_2 \le N$.
We denote by $V_{[ \lambda ]} = V_{[\lambda], \Orth_N}$ 
and $S_{[ \lambda ]} = S_{[\lambda],\Orth_N}$ 
the irreducible representation and its character of $\Orth_N$ corresponding 
to a partition $\lambda \in \Par(\Orth_N)$.
Let $V = \Comp^N$ be the defining representation of $\Orth_N$, and
denote by $H_r$ (resp. $E_r$) 
the character of $\Orth_N$ on $S^r(V)$ (resp. $\twedge^r(V)$).
Then $\twedge^r(V)$ is the irreducible representation corresponding to the one-column partition $(1^r)$, 
while $S^r(V)$ is not irreducible if $r \ge 2$ and the quotient $S^r(V)/S^{r-2}(V)$ 
is the irreducible representation corresponding to the one-row partition $(r)$.
Let $R(\Orth_N)$ be the representation ring of $\Orth_N$ 

For an arbitrary partition $\lambda$, 
we define the corresponding \emph{orthogonal Schur function} $s_{[\lambda]} \in \Lambda$ by
putting
\begin{equation}
\label{eq:O-Schur}
s_{[\lambda]}
 =
\det \left(
 h_{\lambda_i-i+j} - h_{\lambda_i-i-j}
\right)_{1 \le i, j \le l(\lambda)}.
\end{equation}
Let $\pi_{\Orth_N} : \Lambda \to R(\Orth_N)$ be the ring homomorphism defined by 
$\pi_{\Orth_N}(h_r) = H_r$ for $r \ge 1$.
The image of an orthogonal Schur function under $\pi_{\Orth_N}$ can be 
expressed as a linear combination of irreducible characters 
(in fact, it is $0$ or an irreducible character up to sign), 
by using the following algorithm (see \cite{King} and \cite{KT}).

\begin{prop}
\label{prop:O-algorithm}
Let $\lambda$ be a partition.
\begin{enumerate}
\item[(1)]
If $\lambda \in \Par(\Orth_N)$, then we have
$$
\pi_{\Orth_N} \left( s_{[ \lambda ]} \right)
 = 
S_{[ \lambda ]}.
$$
\item[(2)]
In general, the image $\pi_{\Orth_N} \left( s_{[\lambda]} \right)$ can be computed as follows.
We put $r = \lambda_1$ and
$$
\alpha = (\lambda'_1, \lambda'_2-1, \dots, \lambda'_r-(r-1)).
$$
\begin{enumerate}
\item[(a)]
If $\alpha$ has an entry larger than or equal to $N+r$, 
then $\pi_{\Orth_N} \left( s_{[\lambda]} \right) = 0$.
\item[(b)]
If $\alpha_i + \alpha_j = N$ for some $i$ and $j$, 
then $\pi_{\Orth_N} \left( s_{[\lambda]} \right) = 0$.
\item[(c)]
Otherwise, suppose that $\alpha_1 > \dots > \alpha_p > N/2 \ge \alpha_{p+1}$ 
and define a sequence $\beta$ by putting
\begin{multline*}
\beta
\\
 =
\begin{cases}
(N-\alpha_1, \dots, N-\alpha_p, \alpha_{p+1}, \dots, \alpha_r)
 &\text{if $p$ is even,} \\
(N-\alpha_1, \dots, N-\alpha_{p+1}, \alpha_{p+2}, \dots, \alpha_r)
 &\text{if $p$ is odd and $\alpha_p + \alpha_{p+1} \ge N+1$,} \\
(N-\alpha_1, \dots, N-\alpha_{p-1}, \alpha_p, \dots, \alpha_r)
 &\text{if $p$ is odd and $\alpha_p + \alpha_{p+1} \le N-1$.}
\end{cases}
\end{multline*}
Let $\gamma$ be the rearrangement of $\beta$ in decreasing order 
and $\sigma \in \Sym_r$ a permutation satisfying $\sigma(\beta) = \gamma$.
If $\mu$ is the partition given by $\gamma = (\mu'_1, \mu'_2-1, \dots, \mu'_r-(r-1))$, 
then $\mu \in \Par(\Orth_N)$ and 
$$
\pi_{\Orth_N} \left( s_{[\lambda]} \right)
 = 
\sgn(\sigma) S_{[\mu]}.
$$
\end{enumerate}
\item[(3)]
In particular,
if $\trans\lambda_1 + \trans\lambda_2 = N+1$ or $\trans\lambda_1 + \trans\lambda_3 = N+2$, 
then we have $\pi_{\Orth_N} \left( s_{[\lambda]} \right) = 0$.
\end{enumerate}
\end{prop}

We have the following relations in the ring of symmetric functions, 
from which we can derive identities involving irreducible characters of $\Orth_N$ 
by applying $\pi_{\Orth_N}$.

\begin{prop}
\label{prop:O-schur-rel}
\begin{enumerate}
\item[(1)]
(Newell, Littlewood)
For any partitions $\mu$ and $\nu$, we have
\begin{equation}
\label{eq:O-LR}
s_{[\mu]} \cdot s_{[\nu]}
 =
\sum_{\lambda \in \Par}
 \left(
  \sum_{\tau, \, \xi, \, \eta \in \Par} \LR^\mu_{\tau,\xi} \LR^\nu_{\tau,\eta} \LR^\lambda_{\xi,\tau}
 \right)
 s_{[\lambda]}.
\end{equation}
\item[(2)]
(Littlewood)
For any partition $\lambda$, we have
\begin{equation}
\label{eq:O-GL}
s_\lambda
 =
\sum_{\mu \in \Par}
 \left(
  \sum_{\kappa \in \mathcal{E}} \LR^\lambda_{\kappa,\mu}
 \right)
 s_{[\mu]},
\end{equation}
where $\mathcal{E}$ is the set of all partitions whose row lengths are all even.
\end{enumerate}
\end{prop}

Finally we review the representation theory of the special orthogonal group $\SOrth_N$.
We associate to a partition $\lambda \in \Par(\Orth_N)$ another partition $\lambda^{\sharp} \in \Par(\Orth_N)$ 
obtained by replacing the first column (of length $\lambda'_1$) by the column of length $N-\lambda'_1$.
Then we have
$$
V_{[\lambda],\Orth_N} \otimes \Comp_{\det} \cong V_{[\lambda^{\sharp}],\Orth_N},
\quad
\Res^{\Orth_N}_{\SOrth_N} V_{[\lambda],\Orth_N}
 \cong
\Res^{\Orth_N}_{\SOrth_N} V_{[\lambda^{\sharp}],\Orth_N},
$$
where $\Comp_{\det}$ is the one-dimensional representation of $\Orth_N$ given by the determinant.
For a partition $\lambda$ of length $\le N/2$, we denote by $V_{[\lambda],\SOrth_N}$ 
and $S_{[\lambda],\SOrth_N}$ the restriction of $V_{[\lambda],\Orth_N}$ 
and $S_{[\lambda],\Orth_N}$ respectively.
If $N = 2n+1$ is odd, then the restriction $V_{[\lambda],\SOrth_{2n+1}}$ remains irreducible 
and $\{ V_{[\lambda],\SOrth_{2n+1}} : \lambda \in \Par, \ l(\lambda) \le n \}$ forms 
a complete set of representatives of isomorphism classes of irreducible representations 
of $\SOrth_{2n+1}$.
The even orthogonal group case is more subtle.
It is known that the irreducible representations of $\SOrth_{2n}$ are parametrized by 
sequences $\omega = (\omega_1, \dots, \omega_{n-1}, \omega_n)$ of integers satisfying 
$\omega_1 \ge \dots \ge \omega_{n-1} \ge |\omega_n|$.
We denote the corresponding irreducible representation by $L_{[\omega],\SOrth_{2n}}$.
If $\lambda$ is a partition of length $<n$, then the restriction $V_{[\lambda],\SOrth_{2n}}$ 
is the irreducible representation $L_{[(\lambda_1, \dots, \lambda_{l(\lambda)}, 0, \dots, 0)],\SOrth_{2n}}$.
If $\lambda$ is a partition of length $n$, then 
$V_{[\lambda],\SOrth_{2n}}$ is not irreducible and 
decomposes into the direct sum of two distinct irreducible representations 
$L_{[\lambda^+],\SOrth_{2n}}$ and $L_{[\lambda^-],\SOrth_{2n}}$ corresponding to 
$\lambda^+ = (\lambda_1, \dots, \lambda_{n-1}, \lambda_n)$ and 
$\lambda^- = (\lambda_1, \dots, \lambda_{n-1}, - \lambda_n)$.

Let $\{ f^B_r(x) \}$ and $\{ f^D_r(x) \}$ be the sequences of Laurent polynomials defined by
$$
f^B_r(x)
 = 
\begin{cases}
 x^{r+1/2} - x^{-r-1/2} &\text{if $r \ge 0$, } \\
 0 &\text{if $r < 0$,}
\end{cases}
\quad
f^D_r(x)
 =
\begin{cases}
 x^r + x^{-r} &\text{if $r > 0$, } \\
 1 &\text{if $r = 0$, } \\
 0 &\text{if $r < 0$.}
\end{cases}
$$
For a sequence $\alpha = (\alpha_1, \dots, \alpha_n)$ of integers 
and a sequence $\vectx = (x_1, \dots, x_n)$ of indeterminates, 
we put
$$
A^B_\alpha(\vectx)
 =
\left( f^B_{\alpha_j+n-j}(x_i) \right)_{1 \le i, j \le n},
\quad
A^D_\alpha(\vectx)
 =
\left( f^D_{\alpha_j+n-j}(x_i) \right)_{1 \le i, j \le n}.
$$
Then the ratios 
$\det A^B_{\alpha}(\vectx)/\det A^B_{\emptyset}(\vectx)$ and 
$\det A^D_{\alpha}(\vectx)/\det A^D_{\emptyset}(\vectx)$ are Laurent polynomials in $x_1, \dots, x_n$.
Then the Weyl character formula is rephrased as follows:

\begin{prop}
\label{prop:SO-Weyl}
Let $\lambda$ be a partition of length $\le n$.
\begin{enumerate}
\item[(1)]
if $X \in \SOrth_{2n+1}$ has the eigenvalues $x_1, \dots, x_n, x_1^{-1}, \dots, x_n^{-1}$ and $1$, 
then we have
\begin{equation}
\label{eq:SOodd-Weyl}
S_{[ \lambda ],\SOrth_{2n+1}}(X)
 =
\frac{ \det A^B_\lambda(\vectx) }
     { \det A^B_\emptyset(\vectx) }.
\end{equation}
\item[(2)]
If $X \in \SOrth_{2n}$ has the eigenvalues $x_1, \dots, x_n, x_1^{-1}, \dots, x_n^{-1}$, 
then we have
\begin{equation}
\label{eq:SOeven-Weyl}
S_{[ \lambda ],\SOrth_{2n}}(X)
 =
\frac{ \det A^D_\lambda (\vectx) }
     { \det A^D_\emptyset (\vectx) }.
\end{equation}
\end{enumerate}
\end{prop}

\section{%
Pieri rules for the classical groups
}

In this section, we give the Pieri rules for $\Symp_{2n}$ and $\Orth_N$, 
which describe the irreducible decomposition of the tensor product
of an irreducible representation with the symmetric power 
of the defining representation.
The proof uses the symplectic and orthogonal Schur functions and 
their specialization algorithms

At the level of symplectic and orthogonal Schur functions, 
we have the following ``universal'' Pieri rule.

\begin{prop}
\label{prop:univPieri}
For a partition $\mu$ and a nonnegative integer $r$, we have
\begin{gather}
\label{eq:Sp-univPieri}
s_{\langle \mu \rangle} \cdot s_{\langle (r) \rangle}
 =
\sum_{\lambda \in \Par} \# \mathcal{M}^\lambda_{\mu,r} s_{\langle \lambda \rangle},
\\
\label{eq:O-univPieri}
s_{[ \mu ]} \cdot s_{[ (r) ]}
 =
\sum_{\lambda \in \Par} \# \mathcal{M}^\lambda_{\mu,r} s_{\langle \lambda \rangle},
\end{gather}
where $\mathcal{M}^\lambda_{\mu,r}$ is given by
\begin{equation}
\label{eq:univPieri-coeff}
\mathcal{M}^\lambda_{\mu,r}
 =
\{ \xi \in \Par : 
 \text{$\mu/\xi$ and $\lambda/\xi$ are both horizontal strips and $|\mu/\xi|+|\lambda/\xi| = r$}
\}.
\end{equation}
\end{prop}

\begin{demo}{Proof}
We use Newell--Littlewood formulas (\ref{eq:Sp-LR}) and (\ref{eq:O-LR}) with $\nu = (r)$.
Since $\LR^{(r)}_{\tau,\eta} = 0$ unless $\tau =(s)$ and $\eta = (r-s)$ for some $0 \le s \le r$, 
the claim follows from the Pieri rule (\ref{eq:GL-Pieri1}).
\end{demo}

By applying the homomorphisms $\pi_{\Symp_{2n}}$ to (\ref{eq:Sp-univPieri}) 
and $\pi_{\Orth_N}$ to (\ref{eq:O-univPieri}), 
and then by using the algorithms given in Propositions~\ref{prop:Sp-algorithm} 
and \ref{prop:O-algorithm}, we can prove the following ``actual'' Pieri rules.
Part (1) of the following theorem was obtained by Sundaram \cite[Theorem~4.1]{Su1}, 
where she used the Berele insertion algorithm to give a combinatorial proof.

\begin{theorem}
\label{thm:Pieri}
\begin{enumerate}
\item[(1)]
Let $\lambda$, $\mu \in \Par(\Symp_{2n})$ and $r$ a nonnegative integer.
Then the multiplicity of the irreducible $\Symp_{2n}$-module $V_{\langle \lambda \rangle}$ 
in the tensor product $V_{\langle \mu \rangle} \otimes V_{\langle (r) \rangle} 
\cong V_{\langle \mu \rangle} \otimes S^r(V)$, where $V$ is the defining representation of $\Symp_{2n}$, 
is equal to the number of partitions $\xi$ satisfying the following two conditions:
\begin{enumerate}
\item[(i)]
$\mu/\xi$ and $\lambda/\xi$ are both horizontal strips.
\item[(ii)]
$|\mu/\xi| + |\lambda/\xi| = r$.
\end{enumerate}
\item[(2)]
Let $\lambda$, $\mu \in \Par(\Orth_N)$ and $r$ a nonnegative integer.
Then the multiplicity of the irreducible $\Orth_N$-module $V_{[ \lambda ]}$ 
in the tensor product $V_{[ \mu ]} \otimes V_{[(r)]}$ is equal to 
the number of partitions $\xi$ satisfying the following three conditions:
\begin{enumerate}
\item[(i)]
$\mu/\xi$ and $\lambda/\xi$ are both horizontal strips.
\item[(ii)]
$|\mu/\xi| + |\lambda/\xi| = r$.
\item[(iii)]
If $\mu'_1 + \mu'_2 = N$, $\lambda'_1 = \mu'_1$ and $\lambda'_2 = \mu'_2$, then 
one of the following holds:
\begin{enumerate}
\item[(iii-1)]
$l(\xi) = l(\mu)$ and $\xi_l \in \{ \mu_l, \lambda_l \}$, where $l = l(\mu) = l(\xi) = l(\lambda)$.
\item[(iii-2)]
$l(\xi) < l(\mu)$ and $\xi'_2 < \mu'_2$.
\end{enumerate}
\end{enumerate}
\end{enumerate}
\end{theorem}

\begin{demo}{Proof}
(1)
By applying $\pi_{\Symp_{2n}}$ to (\ref{eq:Sp-univPieri}) 
and using Proposition~\ref{prop:Sp-algorithm} (1), we have
$$
S_{\langle \mu \rangle} \cdot S_{\langle r \rangle}
 =
\sum_{l(\lambda) \le n}
 \# \mathcal{M}^\lambda_{\mu,r} S_{\langle \lambda \rangle}
 +
\sum_{l(\lambda) \ge n+1}
 \# \mathcal{M}^\lambda_{\mu,r}
 \pi_{\Symp_{2n}} ( s_{\langle \lambda \rangle} ).
$$
If $\mathcal{M}^\lambda_{\mu,r} \neq \emptyset$, then
$\lambda$ is obtained from a subdiagram of $\mu$ by adding a horizontal strip, 
so we have $l(\lambda) \le l(\mu)+1 \le n+1$.
If $l(\lambda) = n+1$, then $\pi_{\Symp_{2n}} ( s_{\langle \lambda \rangle} ) = 0$ 
by Proposition~\ref{prop:Sp-algorithm} (3).
Hence we see that
$$
S_{\langle \mu \rangle} \cdot S_{\langle r \rangle}
 =
\sum_{l(\lambda) \le n} \# \mathcal{M}^\lambda_{\mu,r} S_{\langle \lambda \rangle},
$$
where $\lambda$ runs over all partitions of length $\le n$.

(2)
By applying $\pi_{\Orth_N}$ to (\ref{eq:O-univPieri}) 
and using Proposition~\ref{prop:O-algorithm} (1) and (3), we have
$$
S_{[ \mu ]} \cdot S_{[ r ]}
 =
\sum_{\lambda \in \Par(\Orth_N)} \# \mathcal{M}^\lambda_{\mu,r} S_{[\lambda]}
 +
\sum_{\rho'_1 + \rho'_2 \ge N+2}
 \# \mathcal{M}^\rho_{\mu,r}
 \pi_{\Orth_N} ( s_{[\rho]} ),
$$
where $\rho$ runs over all partitions satisfying $\rho'_1 + \rho'_2 \ge N+2$.

Suppose that a partition $\rho$ satisfies $\rho'_1 + \rho'_2 \ge N+2$ 
and $\mathcal{M}^\rho_{\mu,r} \neq \emptyset$.
Then $\rho$ is obtained from a subdiagram of $\mu \in \Par(\Orth_N)$ by adding a horizontal strip, 
so we have
$$
\mu'_1 + \mu'_2 = N,
 \quad 
\rho'_1 = \mu'_1+1,
 \quad
\rho'_2 = \mu'_2+1
$$
and
$$
\eta'_1 = \mu'_1,
 \quad
\eta'_2 = \mu'_2
\quad\text{for all $\eta \in \mathcal{M}^\rho_{\mu,r}$.}
$$
If $\rho'_2 = \rho'_3$, then $\pi_{\Orth_N} ( s_{[\rho]} ) = 0$ by 
Proposition~\ref{prop:O-algorithm} (3).
Also it follows from Proposition~\ref{prop:O-algorithm} (2) that, 
if $\pi_{\Orth_N} ( s_{[\rho]} ) \neq 0$, then we have
$$
\pi_{\Orth_N} (s_{[\rho]}) = - S_{[\sigma]},
$$
where the partition $\sigma \in \Par(\Orth_N)$ is given by
$$
\sigma' = (\rho'_1 - 1, \rho'_2 - 1, \rho'_3, \dots).
$$
Therefore we have
$$
S_{[ \mu ]} \cdot S_{[ r ]}
 =
\sum_{\lambda \in \Par(\Orth_N)} \# \mathcal{M}^\lambda_{\mu,r} S_{[\lambda]}
 -
\delta_{\mu'_1+\mu'_2,N}
\sum_{\sigma \in \Par(\Orth_N), \sigma'_1 = \mu'_1, \sigma'_2 = \mu'_2} 
 \# \mathcal{M}^{\tilde{\sigma}}_{\mu,r} S_{[\sigma]},
$$
where the second summation is taken over all $\sigma \in \Par(\Orth_N)$ satisfying 
$\sigma'_1 = \mu'_1$ and $\sigma'_2 = \mu'_2$, 
and $\tilde{\sigma}$ is given by 
$\left( \tilde{\sigma} \right)' = (\sigma'_1+1, \sigma'_2+1, \sigma'_3, \dots)$.

We fix two partitions $\lambda$ and $\mu$ such that 
$\mu'_1 + \mu'_2 = N$, $\lambda'_1 = \mu'_1$ and $\lambda'_2 = \mu'_2$, 
and put
$$
\overline{\mathcal{M}}^\lambda_{\mu,r}(N)
 =
\{ \xi \in \mathcal{M}^\lambda_{\mu,r} :
 \text{$\xi$ satisfies the condition (iii) in Theorem~\ref{thm:Pieri} (2)}
\}.
$$
Then we shall show that
$$
\# \mathcal{M}^\lambda_{\mu,r} - \# \mathcal{M}^{\tilde\lambda}_{\mu,r}
 =
\# \overline{\mathcal{M}}^\lambda_{\mu,r}(N).
$$
Let $\phi: \mathcal{M}^{\tilde\lambda}_{\mu,r} \to \mathcal{M}^\lambda_{\mu,r}$ be the map 
defined by
$$
\phi(\eta) = (\eta_1, \dots, \eta_{l(\eta)-1}, \eta_{l(\eta)}-1)
\quad\text{for $\eta \in \mathcal{M}^{\tilde\lambda}_{\mu,r}$.}
$$
Since $\phi$ is injective, it is enough to show that
$$
\mathcal{M}^\lambda_{\mu,r}
 =
\overline{\mathcal{M}}^\lambda_{\mu,r}(N)
\sqcup
\phi \left( \mathcal{M}^{\tilde\lambda}_{\mu,r} \right).
$$

First we consider the case where $\mu'_1 > \mu'_2$.
If $\eta \in \mathcal{M}^{\tilde\lambda}_{\mu,r}$ and $\xi = \phi(\eta)$, 
then $\xi'_1 < \eta'_1 = \mu'_1$ and $\xi'_2 = \eta'_2 = \mu'_2$, 
so $\xi \not\in \overline{\mathcal{M}}^\lambda_{\mu,r}(N)$.
If $\xi \in \mathcal{M}^\lambda_{\mu,r}$ satisfies $\xi'_1 = \mu'_1$, 
then it follows from $\xi \subset \mu$ that $\xi_{l(\xi)} = 1 = \mu'_1 = \lambda'_1$.
Hence, if $\xi \in \mathcal{M}^\lambda_{\mu,r} \setminus \overline{\mathcal{M}}^\lambda_{\mu,r}(N)$, then 
we have $\xi'_1 < \mu'_1$ and $\xi'_2 = \mu'_2$, 
so $\xi \in \phi \left( \mathcal{M}^{\tilde\lambda}_{\mu,r} \right)$.

Next we consider the case where $\mu'_1 = \mu'_2$.
In this case, $N = 2n$ is even and $\mu_n \ge 2$, $\lambda_2 \ge 2$.
If $\eta \in \mathcal{M}^{\tilde\lambda}_{\mu,r}$ and $\xi = \phi(\eta)$, 
then $\eta_n \ge 2$ (since $\eta'_1 = \eta'_2$) and $\eta_n \le \min \{ \mu_n, \lambda_n \}$,
so $0 < \xi_n < \min \{ \mu_n, \lambda_n \}$ and $\xi \not\in \overline{\mathcal{M}}^\lambda_{\mu,r}(N)$.
If $\xi \in \mathcal{M}^\lambda_{\mu,r} \setminus \overline{\mathcal{M}}^\lambda_{\mu,r}(N)$, 
then we have $\xi'_1 = \mu'_1 = n$ and $\xi_n \not\in \{ \mu_n, \lambda_n \}$, 
so $0 < \xi_n < \min \{ \mu_n, \lambda_n \}$ and 
$\xi \in \phi \left( \mathcal{M}^{\tilde\lambda}_{\mu,r} \right)$.

This completes the proof of (2).
\end{demo}

In the orthogonal group case, the symmetric power $S^r(V)$ of the defining representation $V$ 
of $\Orth_N$ is decomposed as follows:
$$
S^r(V) \cong \bigoplus_{s=0}^{\lfloor r/2 \rfloor} V_{[(r-2s)]},
$$
where $\lfloor r/2 \rfloor$ is the largest integer not exceeding $r/2$.
Hence we have

\begin{corollary}
\label{cor:O-Pieri}
Let $\lambda$, $\mu \in \Par(\Orth_N)$ and $r$ a nonnegative integer.
Then the multiplicity of $V_{[ \lambda ]}$ 
in the tensor product $V_{[ \mu ]} \otimes S^r(V)$, 
where $V$ is the defining representation of $\Orth_N$, is equal to 
the number of partitions $\xi$ satisfying the following three conditions:
\begin{enumerate}
\item[(i)]
$\mu/\xi$ and $\lambda/\xi$ are both horizontal strips.
\item[(ii)]
$|\mu/\xi| + |\lambda/\xi| = r-2s$ for some integer $0 \le s \le r/2$.
\item[(iii)]
If $\mu'_1 + \mu'_2 = N$, $\lambda'_1 = \mu'_1$ and $\lambda'_2 = \mu'_2$, then 
one of the following holds:
\begin{enumerate}
\item[(iii-1)]
$l(\xi) = l(\mu)$ and $\xi_l \in \{ \mu_l, \lambda_l \}$, where $l = l(\mu) = l(\xi) = l(\lambda)$.
\item[(iii-2)]
$l(\xi) < l(\mu)$ and $\xi'_2 < \mu'_2$.
\end{enumerate}
\end{enumerate}
\end{corollary}

Also we have the following Pieri rules for the special orthogonal groups.
Part (1) of the following corollary was given in \cite[Theorem~5.3]{Su2}.

\begin{corollary}
\label{cor:SO-Pieri}
\begin{enumerate}
\item[(1)]
Let $\lambda$ and $\mu$ be partitions of length $\le n$ and $r$ a nonnegative integer.
Then the multiplicity of the irreducible $\SOrth_{2n+1}$-module $V_{[\lambda]}$ 
in $V_{[\mu],\SOrth_{2n+1}} \otimes V_{[ (r) ],\SOrth_{2n+1}}$ 
is equal to the number of partitions $\xi$ satisfying the following three conditions:
\begin{enumerate}
\item[(i)]
$\mu/\xi$ and $\lambda/\xi$ are both horizontal strips.
\item[(ii)]
$|\mu/\xi| + |\lambda/\xi| = r$ or $r-1$.
\item[(iii)]
If $|\mu/\xi| + |\lambda/\xi| = r-1$, then $l(\xi) = l(\mu) = n$.
\end{enumerate}
\item[(2)]
For two partitions $\lambda$ and $\mu$ of length $\le n$ and a nonnegative integer $r$, 
we denote by $\overline{\mathcal{M}}^\lambda_{\mu,r}(N)$ 
the set of partitions $\xi$ satisfying the conditions in Theorem~\ref{thm:Pieri} (2).
Then we have
$$
S_{[\mu],\SOrth_{2n}} \cdot S_{[(r)],\SOrth_{2n}}
\\
 =
\sum_{l(\lambda) \le n} m(\lambda,\mu,n) \# \overline{\mathcal{M}}^\lambda_{\mu,r}(2n) S_{[\lambda],\SOrth_{2n}},
$$
where $m(\lambda,\mu,n)$ is given by
\begin{equation}
\label{eq:mult-SO}
m(\lambda,\mu,n)
 =
\begin{cases}
 2 &\text{if $l(\mu) = n$ and $l(\lambda) < n$, }
\\
 1 &\text{otherwise.}
\end{cases}
\end{equation}
\end{enumerate}
\end{corollary}

Note that, if $l(\mu) < n$ or $l(\lambda) < n$, then $\overline{\mathcal{M}}^\lambda_{\mu,r}(2n) 
= \mathcal{M}^\lambda_{\mu,r}$, where $\mathcal{M}^\lambda_{\mu,r}$ is defined by 
(\ref{eq:univPieri-coeff}).

\begin{demo}{Proof}
(1)
If $l(\mu) \le n$, then $\mu'_1 + \mu'_2 \le 2n < 2n+1$, 
so it follows from Theorem~\ref{thm:Pieri} (2) that
$$
S_{[\mu],\Orth_{2n+1}} \cdot S_{[(r)],\Orth_{2n+1}}
 = 
\sum_{\lambda \in \Par(\Orth_{2n+1})} \# \mathcal{M}^\lambda_{\mu,r} S_{[\lambda],\Orth_{2n+1}}.
$$
If $\rho \in \Par(\Orth_{2n+1})$ satisfies $l(\rho) > n$ and $\mathcal{M}^\rho_{\mu,r} \neq \emptyset$, 
then we have $l(\mu) = n$, $l(\rho) = n+1$ and $\rho_{n+1} = 1$.
In this case, $\rho^{\sharp} = (\rho_1, \dots, \rho_n)$ 
and the map
$$
\varphi : \mathcal{M}^\rho_{\mu,r} \ni \eta \mapsto \eta \in \mathcal{M}^{\rho^{\sharp}}_{\mu,r-1}
$$
is injective and $\varphi\left( \mathcal{M}^\rho_{\mu,r} \right)
 = \{ \xi \in \mathcal{M}^{\rho^{\sharp}}_{\mu,r-1} : l(\xi) = n \}$.
Hence we see that, if $l(\mu) = l(\lambda) = n$, 
then the coefficient of $S_{[\lambda],\SOrth_{2n+1}}$ in $S_{[\mu],\SOrth_{2n+1}} \cdot S_{[r],\SOrth_{2n+1}}$ 
is equal to 
$\# \mathcal{M}^\lambda_{\mu,r} + \# \{ \xi \in \mathcal{M}^\lambda_{\mu,r-1} : l(\xi) = n \}$.

(2)
If $\rho \in \Par(\Orth_{2n})$ satisfies $l(\rho) > n$ and $\overline{\mathcal{M}}^\rho_{\mu,r} \neq \emptyset$, 
then we have $l(\mu) = n$, $l(\rho) = n+1$ and $\rho_n = \rho_{n+1} = 1$.
In this case $\rho^{\sharp} = (\rho_1, \dots, \rho_{n-1})$ and
the map
$$
\varphi : 
\overline{\mathcal{M}}^\rho_{\mu,r}(2n) = \mathcal{M}^\rho_{\mu,r}
 \ni \eta
 \mapsto (\eta_1, \dots, \eta_{n-1}) \in 
\mathcal{M}^{\rho^{\sharp}}_{\mu,r} = \overline{\mathcal{M}}^{\rho^{\sharp}}_{\mu,r}(2n)
$$
is bijective.
The proof follows from this observation.
\end{demo}

\begin{remark}
For the \emph{connected} classical groups $\Symp_{2n}$ and $\SOrth_N$, 
the Pieri rules given in Theorem~\ref{thm:Pieri} (2) and Corollary~\ref{cor:SO-Pieri} 
can be derived by applying the generalized Littlewood--Richardson rule \cite{N}, 
which is obtained from the theory of crystal bases.
By using the generalized Littlewood--Richardson rule, we can show that, 
for given two integer sequences $\lambda = (\lambda_1, \dots, \lambda_n)$ and $\mu = (\mu_1, \dots, \mu_n)$ 
with $\lambda_1 \ge \dots \ge \lambda_{n-1} \ge |\lambda_n|$ 
and $\mu_1 \ge \dots \ge \mu_{n-1} \ge |\mu_n|$, 
the multiplicity of the irreducible $\SOrth_{2n}$-module $L_{[\lambda],\SOrth_{2n}}$ 
in $L_{[\mu],\SOrth_{2n}} \otimes L_{[(r,0,\dots,0)],\SOrth_{2n}}$ is equal to the number of 
integer sequences $\xi$ satisfying the following four conditions:
\begin{enumerate}
\item[(i)]
$\xi_1 \ge \xi_2 \ge \dots \ge \xi_{n-1} \ge |\xi_n|$.
\item[(ii)]
$\mu_1 \ge \xi_1 \ge \mu_2 \ge \xi_2 \ge \dots \ge \xi_{n-1} \ge \mu_n \ge \xi_n$,
and 
$\lambda_1 \ge \xi_1 \ge \lambda_2 \ge \xi_2 \ge \dots \ge \xi_{n-1} \ge \lambda_n \ge \xi_n$.
\item[(iii)]
$\sum_{i=1}^n (\mu_i - \xi_i) + \sum_{i=1}^n (\lambda_i - \xi_i) = r$.
\item[(iv)]
$\xi_n \in \{ \mu_n, \lambda_n \}$.
\end{enumerate}
\end{remark}

In Theorem~\ref{thm:Pieri} and Corollary~\ref{cor:SO-Pieri}, we specialize $r = 1$ to obtain 
the following decomposition of the tensor product with the defining representation.

\begin{corollary}
\label{cor:Pieri-standard}
\begin{enumerate}
\item[(1)]
If $\lambda \in \Par(\Symp_{2n})$ and $V$ is the defining representation of $\Symp_{2n}$, 
then we have
$$
V_{\langle \mu \rangle,\Symp_{2n}} \otimes V
 \cong
\bigoplus_\lambda V_{\langle \lambda \rangle,\Symp_{2n}},
$$
where the direct sum is taken over all $\lambda \in \Par(\Symp_{2n})$ such that 
the diagram of $\lambda$ is obtained from that of $\mu$ by adding or removing one cell.
\item[(2)]
If $\lambda \in \Par(\Orth_N)$ and $V$ is the defining representation of $\Orth_N$, 
then we have
$$
V_{[ \mu ],\Orth_N} \otimes V
 \cong
\bigoplus_\lambda V_{[ \lambda ],\Orth_N},
$$
where the direct sum is taken over all $\lambda \in \Par(\Orth_N)$ such that 
the diagram of $\lambda$ is obtained from that of $\mu$ by adding or removing one cell.
\item[(3)]
If $\lambda$ is a partition of length $\le n$ and $V$ is the defining representation of $\SOrth_{2n+1}$, 
then we have
$$
V_{[ \mu ], \SOrth_{2n+1}} \otimes V
 \cong
\bigoplus_\lambda V_{[ \lambda ], \SOrth_{2n+1}},
$$
where $\lambda$ runs over all partitions of length $\le n$ satisfying one of the following three conditions:
\begin{enumerate}
\item[(i)]
$\lambda \supset \mu$ and $|\lambda| = |\mu| + 1$.
\item[(ii)]
$\lambda \subset \mu$ and $|\lambda| = |\mu| - 1$.
\item[(iii)]
$\lambda = \mu$ and $l(\mu) = n$.
\end{enumerate}
\item[(4)]
If $\lambda$ is a partition of length $\le n$ and $V$ is the defining representation of $\SOrth_{2n}$, 
then we have
$$
V_{[ \mu ], \SOrth_{2n}} \otimes V
 \cong
\bigoplus_\lambda V_{[ \lambda ], \SOrth_{2n}}^{\oplus m(\lambda,\mu,n)},
$$
where $\lambda$ runs over all partitions of length $\le n$ and 
$m(\lambda,\mu,n)$ is given by (\ref{eq:mult-SO}).
\end{enumerate}
\end{corollary}
\section{%
Dual Pieri rules for classical groups
}

In this section, we give the dual Pieri rules for the classical groups, 
which describe the irreducible decomposition of the tensor product of 
an irreducible representation with the exterior power of the defining representation.
We use the Weyl character formulas to obtain the following dual Pieri rules.

\begin{theorem}
\label{thm:dualPieri}
\begin{enumerate}
\item[(1)]
Let $\mu$, $\lambda$ be partitions of length $\le n$ and $r$ an integer with $0 \le r \le 2n$.
The multiplicity of the irreducible $\Symp_{2n}$-module $V_{\langle \lambda \rangle}$ 
in the tensor product $V_{\langle \mu \rangle,\Symp_{2n}} \otimes \twedge^r(V)$, 
where $V$ is the defining representation of $\Symp_{2n}$, 
is equal to 
the number of partitions $\xi$ satisfying the following three conditions:
\begin{enumerate}
\item[(i)]
$l(\xi) \le n$.
\item[(ii)]
$\xi/\mu$ and $\xi/\lambda$ are both vertical strips.
\item[(iii)]
$|\xi/\mu| + |\xi/\lambda| = r$.
\end{enumerate}
\item[(2)]
Let $\mu$, $\lambda$ be partitions of length $\le n$ and $r$ an integer with $0 \le r \le 2n$.
The multiplicity of the irreducible $\SOrth_{2n+1}$ module $V_{[ \lambda ],\SOrth_{2n+1}}$ 
in the tensor product $V_{[ \mu ],\SOrth_{2n+1}} \otimes \twedge^r(V)$, 
where $V$ is the defining representation of $\SOrth_{2n+1}$, 
is equal to 
the number of partitions $\xi$ satisfying the following four conditions:
\begin{enumerate}
\item[(i)]
$l(\xi) \le n$.
\item[(ii)]
$\xi/\mu$ and $\xi/\lambda$ are both vertical strips.
\item[(iii)]
$|\xi/\mu| + |\xi/\lambda| = r$ or $r-1$.
\item[(iv)]
if $l(\mu) < n$, then one of the following holds:
\begin{enumerate}
\item[(iv-1)]
$|\xi/\mu| + |\xi/\lambda| = r$ and $l(\mu) = l(\xi) > l(\lambda)$.
\item[(iv-2)]
$|\xi/\mu| + |\xi/\lambda| = r$ and $l(\xi) = l(\lambda)$.
\item[(iv-3)]
$|\xi/\mu| + |\xi/\lambda| = r-1$ and $l(\xi) = n$.
\end{enumerate}
\end{enumerate}
\item[(3)]
For two partitions $\mu$ and $\lambda$ of length $\le n$ and an integer $r$ with $0 \le r \le 2n$, 
let $\mathcal{K}^\lambda_{\mu,r}(n)$ be the number of partitions $\xi$ satisfying the following four conditions:
\begin{enumerate}
\item[(i)]
$l(\xi) \le n$.
\item[(ii)]
$\xi/\mu$ and $\xi/\lambda$ are both vertical strips.
\item[(iii)]
$|\xi/\mu| + |\xi/\lambda| = r$.
\item[(iv)]
$l(\xi) \in \{ n, l(\mu), l(\lambda) \}$.
\end{enumerate}
Then we have
$$
S_{[\mu],\SOrth_{2n}} \cdot E_{r,\SOrth_{2n}}
 =
\sum_\lambda m(\lambda,\mu,n) \# \mathcal{K}^\lambda_{\mu,r}(n) S_{[\lambda],\SOrth_{2n}},
$$
where $E_{r,\SOrth_{2n}}$ is the character of the exterior power $\twedge^r(V)$ 
of the defining representation $V$ of $\SOrth_{2n}$, 
$\lambda$ runs over all partitions of length $\le n$ and
$m(\lambda,\mu,n)$ is given by (\ref{eq:mult-SO}).
\end{enumerate}
\end{theorem}

\begin{remark}
Part (1) of Theorem~\ref{thm:dualPieri} was given in \cite[Theorem~4.4]{Su1}.
For the special orthogonal groups, 
Sundaram \cite[Theorem~5.4]{Su2} and Weyman \cite[Theorems~$B_n$ and $D_n$]{W} 
gave similar dual Pieri formulas.
It is also possible to apply the generalized Littlewood--Richardson rule \cite{N} 
to obtain dual Pieri rules, but the resulting formulas look more complicated than the formulas 
presented in Theorem~\ref{thm:dualPieri}.
\end{remark}

\begin{demo}{Proof of Theorem~\ref{thm:dualPieri}}
For two partitions $\mu$ and $\lambda$ and nonnegative integers $r$ and $n$, we put
$$
\mathcal{N}^\lambda_{\mu,r}(n)
 =
\# \{ \xi \in \Par : \text{$l(\xi) \le n$, $\xi/\mu$ and $\xi/\lambda$ are vertical strips and
 $|\xi/\mu| + |\xi/\lambda| = r$} \}.
$$

(1)
Suppose that $X \in \Symp_{2n}$ has the eigenvalues $x_1, \dots, x_n, x_1^{-1}, \dots, x_n^{-1}$.
We shall show that
$$
S_{\langle \mu \rangle}(X) \cdot \sum_{r=0}^{2n} E_r(X) t^r
 =
\sum_{l(\lambda) \le n} 
 \left( \sum_{r=0}^{2n} \# \mathcal{N}^\lambda_{\mu,r} t^r \right) 
 S_{\langle \lambda \rangle}(X).
$$

Since the generating function of $E_r(X)$, $0 \le r \le 2n$, is given by
$$
\sum_{r=0}^{2n} E_r(X) t^r
 =
\prod_{i=1}^n (1 + x_i t)(1 + x_i^{-1} t),
$$
it follows from the Weyl character formula (\ref{eq:Sp-Weyl}) that
$$
S_{\langle \mu \rangle}(X) \cdot \sum_{r=0}^N E_r(X) t^r
 =
\frac{ 1 }
     { \det A^C_\emptyset (\vectx) }
\det \left( (1+x_it)(1+x_i^{-1}t) f^C_{\mu_j+n-j}(x_i) \right)_{1 \le i, j \le n}.
$$
By using the relation
$$
(1+x t)(1+x^{-1}t) f^C_r(x)
 =
t f^C_{r-1}(x) + (1+t^2) f^C_r(x) + t f^C_{r+1}(x),
$$
we see that
$$
S_{\langle \mu \rangle}(X) \cdot \sum_{r=0}^{2n} E_r(X) t^r
 =
\frac{ 1 }
     { \det A^C_\emptyset (\vectx) }
\sum_{ \ep, \, \delta \in \{ 0, 1 \}^n }
 \det A^C_{\mu + \ep - \delta}(\vectx) t^{|\ep|+|\delta|},
$$
where $|\ep| = \sum_{i=1}^n \ep_i$ and $|\delta| = \sum_{i=1}^n \delta_i$.
We divide the summation into three parts.
We put
\begin{align*}
\mathcal{A}
 &=
\{ (\ep, \delta) \in \{ 0, 1 \}^n \times \{ 0, 1 \}^n : 
\text{$\mu+\ep$ and $\mu+\ep-\delta$ are partitions} \},
\\
\mathcal{B}
 &=
\{ (\ep, \delta) \in \{ 0, 1 \}^n \times \{ 0, 1 \}^n : 
\text{$\mu+\ep$ is a partition but $\mu+\ep-\delta$ is not a partition} \},
\\
\mathcal{C}
 &=
\{ (\ep, \delta) \in \{ 0, 1 \}^n \times \{ 0, 1 \}^n : 
\text{$\mu+\ep-\delta$ is not a partition} \}.
\end{align*}

If $(\ep, \delta) \in \mathcal{B}$ and $\xi = \mu+\ep$, $\rho = \mu+\ep-\delta$, 
then there exists an index $j$ such that
$$
\xi_j = \xi_{j+1}, \quad \delta_j = 1, \quad \delta_{j+1} = 0.
$$
In this case, $\rho_j + 1 = \rho_{j+1}$ and 
the $j$th and $(j+1)$st columns of the matrix $A^C_{\mu+\ep-\delta}(\vectx)$ are identical, 
so we have $\det A^C_{\mu+\ep-\delta}(\vectx) = 0$.

It remains to show that
$$
\sum_{(\ep, \delta) \in \mathcal{C}}
 \det A^C_{\mu + \ep - \delta}(\vectx) t^{|\ep|+|\delta|}
 =
0.
$$
If $\mu+\ep$ is not a partition, then
there exists an index $j$ such that $\mu_j = \mu_{j+1}$, $\ep_j = 0$ and $\ep_{j+1} = 1$.
For $j=1, \dots, n-1$, we put
$$
\mathcal{C}_j
 =
\{ (\ep, \delta) \in \mathcal{C} : 
\text{$j$ is the smallest index such that $\mu_j = \mu_{j+1}$, $\ep_j = 0$ and $\ep_{j+1} = 1$} \}
$$
and
$$
\mathcal{C}_{j,1}
 =
\{ (\ep, \delta) \in \mathcal{C}_j : \delta_j = \delta_{j+1} \},
\quad
\mathcal{C}_{j,2}
 =
\{ (\ep, \delta) \in \mathcal{C}_j : \delta_j \neq \delta_{j+1} \}.
$$
Then we have $\mathcal{C} = \bigsqcup_{j=1}^{n-1} \mathcal{C}_{j,1} \sqcup \bigsqcup_{j=1}^{n-1} \mathcal{C}_{j,2}$.
If $(\ep, \delta) \in \mathcal{C}_{j,1}$ and $\rho = \mu+\ep-\delta$, 
then we have $\rho_j + 1 = \rho_{j+1}$, so $\det A^C_{\mu+\ep-\delta}(\vectx) = 0$. 
In order to show the summation over $\mathcal{C}_{j,2}$ vanishes, 
we introduce an involution on $\mathcal{C}_{j,2}$.
To a pair $(\ep, \delta) \in \mathcal{C}_{j,2}$, we associate another pair $(\ep, \delta^*)$, 
where $\delta^*$ is given by
$$
\delta^* = (\delta_1, \dots, \delta_{j-1}, \delta_{j+1}, \delta_j, \delta_{j+2}, \dots, \delta_n).
$$
If $\rho = \mu+\ep-\delta$ and $\rho^* = \mu+\ep-\delta^*$, then we have 
$\rho_j = \mu_j-1$, $\rho_{j+1} = \mu_j +1$, $\rho^*_j = \mu_j$ and $\rho^*_{j+1} = \mu_j$,
so we see that 
$\det A^C_{\rho^*}(\vectx) = - \det A^C_\rho(\vectx)$.
Hence we have
$$
\sum_{(\ep, \delta) \in \mathcal{C}_{j,2}}
 \det A^C_{\mu + \ep - \delta}(\vectx) t^{|\ep|+|\delta|}
 =
0.
$$
This completes the proof of (1).

(2)
Let $X \in \SOrth_{2n+1}$ have the eigenvalues $x_1, \dots, x_n, x_1^{-1}, \dots, x_n^{-1}, 1$.
By using the Weyl character formula (\ref{eq:SOodd-Weyl}) and
$$
\sum_{r=0}^{2n+1} E_r(X) t^r
 =
(1+t) \prod_{i=1}^n (1 + x_i t)(1 + x_i^{-1} t),
$$
we have
\begin{multline*}
S_{[\mu]}(X) \cdot \sum_{r=0}^{2n+1} E_r(X) t^r
\\
=
(1+t)
 \cdot
\frac{ 1 }
     { \det A^B_\emptyset(\vectx) }
\det \left(
 (1+x_it)(1+x_i^{-1}t) f^B_{\mu_j+n-j}(x_i)
\right)_{1 \le i, j \le n}.
\end{multline*}
We note that
$$
(1+xt)(1+x^{-1}t) f^B_r(x)
 =
\begin{cases}
 (1-t+t^2) f^B_0(x) + t f^B_1(x) &\text{if $r = 0$,} \\
 t f^B_{r-1}(x) + (1+t^2) f^B_r(x) + t f^B_r(x) &\text{if $r \ge 1$.}
\end{cases}
$$

First we consider the case where $l(\mu) = n$. 
In this case, by the same argument as in the proof of (1), 
we have
$$
S_{[\mu]}(X)
\cdot
\prod_{i=1}^n (1 + x_i t)(1 + x_i^{-1} t)
 =
\sum_{l(\lambda) \le n}
 \left( \sum_{r=0}^{2n} \# \mathcal{N}^\lambda_{\mu,r}(n) t^r \right) S_{[\lambda]}(X),
$$
where $\lambda$ runs over all partitions of length $\le n$. 
Hence we have
$$
S_{[\mu]}(X) \cdot E_r(X)
 =
\sum_{l(\lambda) \le n}
 \left( \# \mathcal{N}^\lambda_{\mu,r}(n) + \# \mathcal{N}^\lambda_{\mu,r-1}(n) \right)
 S_{[\lambda]}(X).
$$

Next we consider the case where $l = l(\mu) < n$.
Let $v_r$ and $w_r$ be the column vectors given by
$$
v_r = \left( f^B_r(x_i) \right)_{1 \le i \le n},
\quad
w_r = \begin{cases}
 (1+t^3) v_0 + (t+t^2) v_1 &\text{if $r=0$,} \\
 t v_{r-1} + (1+t^2) v_r + t v_{r+1} &\text{if $r \ge 1$.}
\end{cases}
$$
Then we have
\begin{multline*}
S_{[\lambda]}(X) \cdot \sum_{r=0}^{2n+1} E_r(X) t^r
\\
 =
\frac{ 1 }
     { \det A^B_\emptyset(\vectx) }
\det \begin{pmatrix} 
 w_{\mu_1+n-1} & \cdots & w_{\mu_l+n-l} & w_{n-l-1} & \cdots & w_1 & w_0
\end{pmatrix}.
\end{multline*}
The transition matrix $T$ from $(v_{n-l}, \dots, v_1, v_0)$ to $(w_{n-l-1}, \dots, w_1, w_0)$ 
is given by
$$
T = \begin{pmatrix}
 t      & 0      & 0      & \cdots & 0      & 0      & 0 \\
 1+t^2  & t      & 0      & \cdots & 0      & 0      & 0 \\
 t      & 1+t^2  & t      & \cdots & 0      & 0      & 0 \\
 \vdots & \vdots & \vdots & \ddots & \vdots & \vdots & \vdots \\
 0      & 0      & 0      & \cdots & t      & 0      & 0 \\
 0      & 0      & 0      & \cdots & 1+t^2  & t      & 0 \\
 0      & 0      & 0      & \cdots & t      & 1+t^2  & t+t^2 \\
 0      & 0      & 0      & \cdots & 0      & t      & 1+t^3
\end{pmatrix}.
$$
By computing the determinants of $(n-l) \times (n-l)$ submatrices of $T$, 
we see that
\begin{multline*}
w_{n-l-1} \wedge w_{n-l-2} \wedge \dots \wedge w_1 \wedge w_0
\\
 =
\sum_{s=0}^{n-l}
 \left( t^{n-l-s} + t^{n-l+s+1} \right) 
 v_{n-l} \wedge \dots \wedge v_{s+1} \wedge v_{s-1} \wedge \dots \wedge v_0.
\end{multline*}
By the same argument as in the proof of (1), 
we have
$$
w_{\mu_1+n-1} \wedge \dots \wedge w_{\mu_l+n-l}
 =
\sum_{l(\rho) \le l}
 \left( \sum_{r=0}^{2n+1} \# \mathcal{N}^\rho_{\mu,r}(l) t^r \right)
 v_{\rho_1+n-1} \wedge \dots \wedge v_{\rho_l+n-l},
$$
where $\rho$ runs over all partitions of length $\le l$. 
Also we note that, if $l(\rho) < l$ and $s \neq l$, then we have
$$
v_{\rho_1+n-1} \wedge \dots \wedge v_{\rho_l+n-l} \wedge
 v_{n-l} \wedge \dots \wedge v_{s+1} \wedge v_{s-1} \wedge \dots \wedge v_0
 =
0.
$$
Hence we have
\begin{align}
&
S_{[\mu]}(X) \cdot \sum_{r=0}^{2n+1} E_r(X) t^r
\notag
\\
&\quad
 =
\sum_{l(\lambda) < l}
 \left( 1 + t^{2(n-l)+1} \right)
 \left( \sum_{r=0}^{2n+1} \# \mathcal{N}^\lambda_{\mu,r}(l) t^r \right)
 S_{[\lambda]}(X)
\notag
\\
&\quad\quad
 +
\sum_{s=0}^{n-l}
 \left( t^{n-l-s} + t^{n-l+s+1} \right) 
\sum_{l(\rho) = l}
 \left( \sum_{r=0}^{2n+1} \# \mathcal{N}^\rho_{\mu,r}(l) t^r \right)
 S_{[\rho \cup (1^{n-l-s})]}(X),
\label{eq:demo-SOodd}
\end{align}
where $\rho \cup (1^{n-l-s}) = (\rho_1, \dots, \rho_l, \underbrace{1, \dots, 1}_{n-l-s})$.
Note that $l(\eta) = l$ for all $\eta \in \mathcal{N}^\lambda_{\mu,r}(l)$ or $\mathcal{N}^\rho_{\mu,r}(l)$.

If $l(\lambda) < l$, then tha maps
\begin{gather*}
\psi_1 : \mathcal{N}^\lambda_{\mu,r}(l) \ni \eta 
\longmapsto \eta \in \mathcal{N}^\lambda_{\mu,r}(n),
\\
\psi_2 : \mathcal{N}^\lambda_{\mu,r}(l) \ni \eta 
\longmapsto \eta \cup (1^{n-l}) \in \mathcal{N}^\lambda_{\mu,r+2(n-l)}(n)
\end{gather*}
are injective and
\begin{gather*}
\psi_1 \left( \mathcal{N}^\lambda_{\mu,r}(l) \right)
 =
\{ \xi \in \mathcal{N}^\lambda_{\mu,r}(n) : 
l(\mu) = l(\xi) > l(\lambda) \},
\\
\psi_2 \left( \mathcal{N}^\lambda_{\mu,r}(l) \right)
 = \{ \xi \in \mathcal{N}^\lambda_{\mu,r+2(n-l)}(n) : l(\xi) = n \}.
\end{gather*}
If $l(\rho) = l$ and $0 \le s \le n-l$, then the maps
\begin{gather*}
\psi_3 : \mathcal{N}^\rho_{\mu,r}(l) \ni \eta 
\longmapsto \eta \cup (1^{n-l-s}) \in \mathcal{N}^{\rho \cup (1^{n-l-s})}_{\mu,r+n-l-s}(n),
\\
\psi_4 : \mathcal{N}^\rho_{\mu,r}(l) \ni \eta 
\longmapsto \eta \cup (1^{n-l}) \in \mathcal{N}^{\rho \cup (1^{n-l-s})}_{\mu,r+n-l+s}(n)
\end{gather*}
are injective and
\begin{gather*}
\psi_3 \left( \mathcal{N}^\rho_{\mu,r}(l) \right)
 =
\{ \xi \in \mathcal{N}^{\rho \cup (1^{n-l-s})}_{\mu,r+n-l-s}(n) : 
l(\xi) = l(\lambda) \},
\\
\psi_4 \left( \mathcal{N}^\rho_{\mu,r}(l) \right)
 = \{ \xi \in \mathcal{N}^{\rho \cup (1^{n-l-s})}_{\mu,r+n-l+s}(n) : l(\xi) = n \}.
\end{gather*}
Combining these observations with (\ref{eq:demo-SOodd}) completes the proof of (2).

(3)
Let $X \in \SOrth_{2n}$ have the eigenvalues $x_1, \dots, x_n, x_1^{-1}, \dots, x_n^{-1}$.
By using the Weyl character formula (\ref{eq:SOeven-Weyl}) and
$$
\sum_{r=0}^{2n+1} E_r(X) t^r
 =
\prod_{i=1}^n (1 + x_i t)(1 + x_i^{-1} t),
$$
we have
$$
S_{[ \mu ]}(X) \cdot \sum_{r=0}^{2n} E_r(X) t^r
 =
\frac{ 1 }
     { \det A^D_\emptyset (\vectx) }
\det \left(
 (1+x_it)(1+x_i^{-1}t) f^D_{\mu_j+n-j}(x_i)
 \right)_{1 \le i, j \le n}
$$
And we have
$$
(1+xt)(1+x^{-1}t) f^D_r(x)
 =
\begin{cases}
 (1+t^2) f^D_0(x) + t f^D_1(x) &\text{if $r = 0$,} \\
 2t f^D_0(x) + (1+t^2) f^D_1(x) + t f^D_2(x) &\text{if $r = 1$, } \\
 t f^D_{r-1}(x) + (1+t^2) f^D_r(x) + t f^D_r(x) &\text{if $r \ge 2$.}
\end{cases}
$$

First we consider the case where $\lambda_n \ge 2$.
In this case, by the same argument as in the proof of (1), 
we have
$$
S_{[\mu]}(X) \cdot E_r(X)
 =
\sum_{l(\lambda) \le n} \# \mathcal{N}^\lambda_{\mu,r}(n) S_{[\lambda]}(X),
$$
where $\lambda$ runs over all partitions of length $\le n$.

Next we consider the case where $\lambda_n = 1$.
Then we have
$$
S_{[\mu]}(X) \cdot \sum_{r=0}^{2n} E_r(X) t^r
 =
\sum_{\ep, \delta \in \{ 0, 1 \}^n}
 2^{c(\ep,\delta)}
 t^{|\ep|+|\delta|}
 S_{[\mu+\ep-\delta]}(X),
$$
where
$$
c(\ep, \delta)
 =
\begin{cases}
 1 &\text{if $\ep_n = 0$ and $\delta_n = 1$,} \\
 0 &\text{otherwise.}
\end{cases}
$$
Now, by the argument similar to that in the proof of (1), 
we obtain the desired result.

In what follows we consider the case where $l(\lambda) < n$.
If $l(\lambda) = n-1$, then by the same argument as in the proof of (1), 
we have
$$
S_{[\mu]}(X) \cdot E_r(X)
 =
\sum_\lambda \# \mathcal{N}^\lambda_{\mu,r}(n) S_{[\lambda]}(X).
$$
So we may assume $l(\lambda) \le n-2$, i.e., $\lambda_{n-1} = \lambda_n = 0$.
Let $v_r$ and $w_r$ be the column vectors given by
$$
v_r = \left( f^D_r(x_i) \right)_{1 \le i \le n},
\quad
w_r = t v_{r-1} + (1+t^2) v_r + t v_{r+1},
$$
where $v_{-1} = 0$.
Then we have
\begin{align*}
&
\det \begin{pmatrix}
 w_{\mu_1+n-1} & \cdots & w_{\mu_{n-2}+2} & w_1 + t v_0 & (1+t^2) v_0 + t v_1
\end{pmatrix}
\\
&\quad
 =
\det \begin{pmatrix}
 w_{\mu_1+n-1} & \cdots & w_{\mu_{n-2}+2} & w_1 & (1+t^2) v_0 + t v_1
\end{pmatrix}
\\
&\quad\quad
 +
\det \begin{pmatrix}
 w_{\mu_1+n-1} & \cdots & w_{\mu_{n-2}+2} & t v_0 & (1+t^2) v_0 + t v_1
\end{pmatrix}
\\
&\quad
 =
\det \begin{pmatrix}
 w_{\mu_1+n-1} & \cdots & w_{\mu_{n-2}+2} & w_1 & w_0
\end{pmatrix}
\\
&\quad\quad
 - t^2
\det \begin{pmatrix}
 w_{\mu_1+n-1} & \cdots & w_{\mu_{n-2}+2} & v_1 & v_0
\end{pmatrix}.
\end{align*}
By using the same argument as in the proof of (1), we see that
$$
S_{[\mu]}(X) E_r(X)
 =
\sum_{l(\lambda) \le n} \# \mathcal{N}^\lambda_{\mu,r}(n) S_{[\lambda]}(X)
 -
\sum_{l(\lambda) \le n-2} \# \mathcal{N}^\lambda_{\mu,r-2}(n-2) S_{[\lambda]},
$$
We put
$$
\overline{\mathcal{N}}^\lambda_{\mu,r}(n)
 =
\{ \xi \in \mathcal{N}^\lambda_{\mu,r}(n) : l(\xi) \in \{ n, l(\mu), l(\lambda) \} \}.
$$
Let $\psi : \mathcal{N}^\lambda_{\mu,r-2}(n-2) \to \mathcal{N}^\lambda_{\mu,r}(n)$ be the map 
given by
$$
\psi(\eta) = \eta \cup (1)
 = (\eta_1, \dots \eta_{l(\eta)}, 1)
\quad
\text{for $\eta \in \mathcal{N}^\lambda_{\mu,r-2}(n-2)$.}
$$
Since $\psi$ is injective, it is enough to show that
$$
\mathcal{N}^\lambda_{\mu,r}(n)
 =
\overline{\mathcal{N}}^\lambda_{\mu,r}(n) \sqcup \psi \left( \mathcal{N}^\lambda_{\mu,r-2}(n-2) \right).
$$

If $\eta \in \mathcal{N}^\mu_{\lambda,r-2}(n-2)$ and $\xi = \psi(\eta)$, 
then $l(\xi) = l(\eta)+1$ and it follows from $l(\eta) \le n-2$ and 
$\mu \subset \xi \supset \lambda$ that $\max \{ l(\mu), l(\lambda) \} < l(\xi) < n$, 
so we have $\xi \not\in \overline{\mathcal{N}}^\lambda_{\mu,r}(n)$.
Conversely, 
if $\xi \in \mathcal{N}^\lambda_{\mu,r}(n) \setminus \overline{\mathcal{N}}^\lambda_{\mu,r}(n)$, 
then we see that $\xi_{l(\xi)} = 1$ and $\xi = \psi(\eta)$ with $\eta = (\xi_1, \dots, \xi_{l(\xi)-1}) \in 
\mathcal{N}^\lambda_{\mu,r-2}(n-2)$.
Hence we have 
$\# \mathcal{N}^\lambda_{\mu,r} - \# \mathcal{N}^\lambda_{\mu,r-2}
 = \# \overline{\mathcal{N}}^\lambda_{\mu,r}$, which completes the proof.
\end{demo}

\section{%
Applications to combinatorics of oscillating tableaux
}

In this section, we apply the Pieri rules obtained in the previous sections 
to derive several equinumeration results between 
down-up/up-down tableaux (generalization of oscillating tableaux) 
and column-strict/row-strict tableaux (generalization of standard tableaux).

\begin{definition}
\label{def:strict-tab}
A filling of the diagram of a partition $\lambda$ with positive integers is called a
\emph{column-strict} (resp. \emph{row-strict}) tableau if it satisfies the following two conditions:
\begin{enumerate}
\item[(i)]
Every row is weakly increasing (resp. strictly increasing),
\item[(ii)]
Every columns is strictly increasing (resp. weakly increasing).
\end{enumerate}
Given a columns-strict or row-strict tableau $T$, 
the \emph{weight} of $T$ is defined to be the sequence $(\alpha_1, \alpha_2, \dots)$, 
where $\alpha_i$ is the number of occurrences of $i$ in $T$.
We denote by $\CSTab(\lambda,\alpha)$ (resp. $\RSTab(\lambda,\alpha)$) 
the set of all column-strict (resp. row-strict) tableaux of shape $\lambda$ and weight $\alpha$.
\end{definition}

A column-strict (resp. row-strict) tableau of shape $\lambda$ 
and weight $\alpha = (\alpha_1, \dots, \alpha_k)$ 
is identified with a sequence
$$
\emptyset = \lambda^{(0)} \subset \lambda^{(1)} \subset \dots 
\subset \lambda^{(k-1)} \subset \lambda^{(k)} = \lambda
$$
of partitions such that $\lambda^{(i)}/\lambda^{(i-1)}$ is a horizontal (resp. vertical) 
$\alpha_i$-strip for each $i$.

\begin{definition}
\label{def:DU-UD-tab}
\begin{enumerate}
\item[(1)]
A sequence $( \lambda^{(i)} )_{i=0}^{2k}$ of partitions is called 
a \emph{down-up tableau} of shape $\lambda$ if it satisfies the following two conditions:
\begin{enumerate}
\item[(i)]
$\lambda^{(0)} = \emptyset$ and $\lambda^{(2k)} = \lambda$.
\item[(ii)]
$\lambda^{(2i-2)} \supset \lambda^{(2i-1)} \subset \lambda^{(2i)}$, and 
$\lambda^{(2i-2)}/\lambda^{(2i-1)}$ and $\lambda^{(2i)}/\lambda^{(2i-1)}$ are both horizontal strip.
\end{enumerate}
\item[(2)]
Dually, a sequence $( \lambda^{(i)} )_{i=0}^{2k}$ of partitions is called 
a \emph{up-down tableau} of shape $\lambda$ if it satisfies the following two conditions:
\begin{enumerate}
\item[(i)]
$\lambda^{(0)} = \emptyset$ and $\lambda^{(2k)} = \lambda$.
\item[(ii)]
$\lambda^{(2i-2)} \subset \lambda^{(2i-1)} \supset \lambda^{(2i)}$, and 
$\lambda^{(2i-1)}/\lambda^{(2i-2)}$ and $\lambda^{(2i-1)}/\lambda^{(2i)}$ are both horizontal strip.
\end{enumerate}
\end{enumerate}
We should remark that the terminologies ``down-up tableau'' and ``up-down tableau'' 
are used for different meanings in \cite{Su1}, \cite{Su2} and other literatures.
\end{definition}

Now we are ready to state our main results.

\begin{theorem}
\label{thm:main2}
Let $\alpha = (\alpha_1, \dots, \alpha_k)$ be a sequence of nonnegative integers.
\begin{enumerate}
\item[(1)]
For nonnegative integers $n$ and $m$, the following two sets are equinumerous:
\begin{enumerate}
\item[(a)]
The set of down-up tableaux $( \lambda^{(i)} )_{i=0}^{2k}$ of shape $(m)$ such that 
$l(\lambda^{(i)}) \le n$ for $0 \le i \le 2k$ and
$|\lambda^{(2i-2)}/\lambda^{(2i-1)}| + |\lambda^{(2i)}/\lambda^{(2i-1)}| = \alpha_i$ 
for $1 \le i \le k$.
\item[(b)]
The set of columns-strict tableaux of weight $\alpha$ whose shape $\lambda$ satisfies 
$l(\lambda) \le 2n$ and $c(\lambda) = m$.
\end{enumerate}
\item[(2)]
For nonnegative integers $n$ and $m$, the following two sets are equinumerous:
\begin{enumerate}
\item[(a)]
The set of up-down tableaux $( \lambda^{(i)} )_{i=0}^{2k}$ of shape $(m)$ such that 
$l(\lambda^{(i)}) \le n$ for $0 \le i \le 2k$ and
$|\lambda^{(2i-1)}/\lambda^{(2i-2)}| + |\lambda^{(2i-1)}/\lambda^{(2i)}| = \alpha_i$ 
for $1 \le i \le k$.
\item[(b)]
The set of row-strict tableaux of weight $\alpha$ whose shape $\lambda$ satisfies 
$l(\lambda) \le 2n$ and $c(\lambda) = m$.
\end{enumerate}
\item[(3)]
For nonnegative integers $N$ and $m$, the following two sets are equinumerous:
\begin{enumerate}
\item[(a)]
The set of down-up tableaux $( \lambda^{(i)} )_{i=0}^{2k}$ of shape $(1^m)$ such that
$\lambda^{(i)} \in \Par(\Orth_N)$ (i.e., $(\lambda^{(i)})'_1 + (\lambda^{(i)})'_2 \le N$) 
for $0 \le i \le 2k$,
$|\lambda^{(2i-2)}/\lambda^{(2i-1)}| + |\lambda^{(2i)}/\lambda^{(2i-1)}| 
\in \{ \alpha_i, \alpha_i-2, \alpha_i - 4, \dots\}$ for $0 \le i \le k$ and 
each triple $(\mu, \xi, \lambda) = (\lambda^{(2i-2)}, \lambda^{(2i-1)}, \allowbreak \lambda^{(2i)})$ satisfies 
the condition (iii) in Theorem~\ref{thm:Pieri} (2).
\item[(b)]
The set of columns-strict tableaux of weight $\alpha$ whose shape $\lambda$ satisfies 
$l(\lambda) \le N$ and $r(\lambda) = m$.
\end{enumerate}
\item[(4)]
For nonnegative integers $n$ and $m$, the following two sets are equinumerous:
\begin{enumerate}
\item[(a)]
The set of up-down tableaux $( \lambda^{(i)} )_{i=0}^{2k}$ of shape $(1^m)$ such that 
$l(\lambda^{(i)}) \le n$ for $0 \le i \le 2k$, 
$|\lambda^{(2i-2)}/\lambda^{(2i-1)}| + |\lambda^{(2i)}/\lambda^{(2i-1)}| 
 \in \{ \alpha_i, \alpha_i - 1 \}$ for $1 \le i \le k$, and 
each triple $(\mu, \xi, \lambda) = (\lambda^{(2i-2)}, \lambda^{(2i-1)}, \lambda^{(2i)})$ satisfies 
the condition (iv) in Theorem~\ref{thm:dualPieri} (2).
\item[(b)]
The set of row-strict tableaux of weight $\alpha$ whose shape $\lambda$ satisfies 
$l(\lambda) \le 2n+1$ and $r(\lambda) = m$ or $2n+1-m$.
\end{enumerate}
\item[(5)]
For nonnegative integers $n$ and $m$, the following two two numbers coincide:
\begin{enumerate}
\item[(a)]
The summation
$$
\sum_U 2^{d(U)}
$$
over over all up-down tableaux $U = ( \lambda^{(i)} )_{i=0}^{2k}$ of shape $(1^m)$ 
such that 
$l(\lambda^{(i)}) \le n$ for $0 \le i \le 2k$, 
$|\lambda^{(2i-2)}/\lambda^{(2i-1)}| + |\lambda^{(2i)}/\lambda^{(2i-1)}| = \alpha_i$ 
for $1 \le i \le k$ 
and 
each triple $(\mu, \xi, \lambda) = (\lambda^{(2i-2)}, \lambda^{(2i-1)}, \lambda^{(2i)})$ satisfies 
the condition (iv) in Theorem~\ref{thm:dualPieri} (3).
Here the statistic $d(U)$ is defined by
$$
d(U) = \# \{ i : l(\lambda^{(2i-2)}) = n, \ l(\lambda^{(2i)}) < n \}.
$$
\item[(b)]
The number of row-strict tableaux of weight $\alpha$ whose shape $\lambda$ satisfies 
$l(\lambda) \le 2n$ and $r(\lambda) = m$ or $2n-m$.
\end{enumerate}
\end{enumerate}
Here $c(\lambda)$ (resp. $r(\lambda)$) is the numbers of columns (resp. rows) of odd length.
\end{theorem}

\begin{demo}{Proof}
We consider the classical group $G$ and its representations $T$ and $W$ listed in the following table:
$$
\begin{array}{c|c|c|c}
    & G & T & W \\
\hline
(1) & \Symp_{2n} & S^\alpha(V) & V_{\langle (m) \rangle} \\
\hline
(2) & \Symp_{2n} & \twedge^\alpha(V) & V_{\langle (m) \rangle} \\
\hline
(3) & \Orth_N & S^\alpha(V) & V_{[(1^m)]} \\
\hline
(4) & \SOrth_{2n+1} & \twedge^\alpha(V) & V_{[(1^m)]} \\
\hline
(5) & \SOrth_{2n} & \twedge^\alpha(V) & V_{[(1^m)]}
\end{array}
$$
Here $S^\alpha(V)$ and $\twedge^\alpha(V)$ is defined by 
$$
S^\alpha(V) = S^{\alpha_1}(V) \otimes \dots \otimes S^{\alpha_k}(V),
\quad
\twedge^\alpha(V) = \twedge^{\alpha_1}(V) \otimes \dots \otimes \twedge^{\alpha_k}(V),
$$
where $V$ is the defining representation of $G$.
We compute the ``multiplicity'' $[T:W]$ of $W$ in $T$ in two ways.
Except for the case where $G = \SOrth_{2n}$ and $W = V_{[n],\SOrth_{2n}}$, 
$W$ is an irreducible representation of $G$, and the multiplicity $[W:T]$ is defined as usual.
If $G = \SOrth_{2n}$ and $W = V_{[n],\SOrth_{2n}}$, 
then the character of $T = \twedge^\alpha (V)$ can be expressed as a linear combination of 
$S_{[\lambda],\SOrth_{2n}}$, $l(\lambda) \le n$, and the multiplicity $[T:W]$ is defined to be the coefficient of 
$S_{[n],\SOrth_{2n}}$ in the character of $T$.

On the one hand, by iteratively using the Pieri rules in Theorems~\ref{thm:Pieri} and \ref{thm:dualPieri},
we see that the multiplicity $[T:W]$ is equal to the number of combinatorial objects given in part (a) 
in each case.

On the other hand, the defining representation of $G$ 
is the restriction of the defining representation of $\GL_N$.
Hence, by using the Pieri rule (\ref{eq:GL-Pieri}) for $\GL_N$, we have the following decomposition 
as $\GL_N$-modules:
$$
S^\alpha(V)
 \cong 
\bigoplus_{l(\lambda) \le N} V_\lambda^{\oplus \# \CSTab(\lambda,\alpha)},
\quad
\twedge^\alpha(V)
 \cong 
\bigoplus_{l(\lambda) \le N} V_\lambda^{\oplus \# \RSTab(\lambda,\alpha)},
$$
where 
$\CSTab(\lambda,\alpha)$ (resp. $\RSTab(\lambda,\alpha)$) denotes the set of 
column-strict (resp. row-strict) tableaux of shape $\lambda$ and weight $\alpha$.

Now by using the following restriction multiplicity formulas (Theorem~\ref{thm:rest}) 
and the relation $\Res^{\Orth_N}_{\SOrth_N} V_{[1^m]} \cong \Res^{\Orth_N}_{\SOrth_N} V_{[(1^{N-m})]}$, 
we can complete the proof of Theorem~\ref{thm:main2}
\end{demo}

\begin{theorem}
\label{thm:rest}
\begin{enumerate}
\item[(1)]
If $\lambda$ is a partition of length $\le 2n$, then the multiplicity of 
$V_{\langle (m) \rangle,\Symp_{2n}}$ in the restriction $\Res^{\GL_{2n}}_{\Symp_{2n}} V_{\lambda,\GL_{2n}}$ 
is given by
$$
\left[
 \Res^{\GL_{2n}}_{\Symp_{2n}} V_\lambda : V_{\langle (m) \rangle}
\right]
 =
\begin{cases}
 1 &\text{if $c(\lambda) = m$,} \\
 0 &\text{otherwise.}
\end{cases}
$$
\item[(2)]
If $\lambda$ is a partition of length $\le N$, then the multiplicity of 
$V_{[(1^m)],\Orth_N}$ in the restriction $\Res^{\GL_N}_{\Orth_N} V_{\lambda,\GL_N}$ 
is given by
$$
\left[
 \Res^{\GL_N}_{\Orth_N} V_\lambda : V_{[(1^m)]}
\right]
 =
\begin{cases}
 1 &\text{if $r(\lambda) = m$,} \\
 0 &\text{otherwise.}
\end{cases}
$$
\end{enumerate}
\end{theorem}

\begin{demo}{Proof}
(1)
It follows from (\ref{eq:Sp-GL}) that the restriction of the character $S_{\lambda,\GL_{2n}}$ 
to $\Symp_{2n}$ can be expressed as
\begin{equation}
\label{eq:Sp-restriction}
\Res^{\GL_{2n}}_{\Symp_{2n}} S_\lambda
 =
\sum_\mu
 \left(
  \sum_{\kappa \in \mathcal{E}'} \LR^\lambda_{\kappa,\mu}
 \right)
 \pi_{\Symp_{2n}}( s_{\langle \mu \rangle} ),
\end{equation}
where $\mu$ runs over all partitions of length $\le 2n$.

Here we use Proposition~\ref{prop:Sp-algorithm} to show 
that $\mu = (m)$ is the only partition of length $\le 2n$ satisfying 
$\pi_{\Symp_{2n}} (s_{\langle \mu \rangle}) = \pm S_{\langle (m) \rangle}$. 
Suppose that $\pi_{\Symp_{2n}} (s_{\langle \mu \rangle}) = \pm S_{\langle (m) \rangle}$.
Let $\alpha = (\mu'_1, \mu'_2 -1, \dots, \mu'_r - r+1)$, 
where $r = \mu_1$.
If $\alpha$ has an entry equal to $n+1$, 
then have $\pi_{\Symp_{2n}}( s_{\langle \mu \rangle} ) = 0$.
We consider the case where $\alpha$ has an entry greater than $n+1$.
In this case, suppose that $\alpha_1 > \dots > \alpha_p > n+1$ and 
put $\beta = (2n+2-\alpha_1, \dots, 2n+2-\alpha_p, \alpha_{p+1}, \dots, \alpha_r)$.
Since $\pi_{\Symp_{2n}}( s_{\langle \mu \rangle} ) = \pm S_{\langle (m) \rangle}$, 
the rearrangement of $\beta$ in decreasing order is equal to $(1,0,-1,\dots,-m+2, -m, \dots, -r+1)$.
However, since $\alpha_1 = \mu'_1 \le 2n$, we have $\beta_1 = 2n+2 - \alpha_1 \ge 2$ 
and this leads a contradiction.
Hence we conclude that $l(\mu) \le n$.
If $l(\mu) \le n$ and $\pi ( s_{\langle \mu \rangle} ) = S_{\langle \mu \rangle}$.
Since the irreducible characters $S_{\langle \nu \rangle}$, $l(\nu) \le n$, are linearly independent, 
we have $\mu = (m)$.

Now it follows from (\ref{eq:Sp-restriction}) that
$$
\left[
 \Res^{\GL_{2n}}_{\Symp_{2n}} V_\lambda : V_{\langle (m) \rangle}
\right]
 =
\sum_{\kappa \in \mathcal{E}'} \LR^\lambda_{\kappa,(m)}.
$$
By the Pieri rule (\ref{eq:GL-Pieri1}) for $\GL_{2n}$, we see that, 
if $\kappa \in \mathcal{E}'$ and $\LR^\lambda_{\kappa,(m)} \neq 0$, then $c(\lambda) = m$, 
and that, if $c(\lambda) = m$, then there is exactly one $\kappa \in \mathcal{E}'$ such that 
$\lambda$ is obtained by adding a horizontal $m$-strip to $\kappa$.
This concludes the proof of (1).

(2)
The proof is similar to that of (1).
We need to prove that 
$\mu = (1^m)$ is the only partition of length $\le N$ satisfying 
$\pi_{\Orth_N} (s_{[ \mu ]}) = \pm S_{[ (1^m) ]}$. 

It is enough to show that $\pi_{\Orth_N} (s_{[ \mu ]}) = \pm S_{[ (1^m)]}$ implies 
$\mu'_1 + \mu'_2 \le N$.
Assume that $\mu'_1 + \mu'_2 > N$ to the contrary.
Let $\alpha = (\mu'_1, \mu'_2 -1, \dots, \mu'_r - r+1)$, 
where $r = \mu_1$.
If $\alpha_i + \alpha_j = N$ for some $i$ and $j$, then $\pi_{\Orth_N}( s_{[\mu]} ) = 0$.
Hence we have $\alpha_1 + \alpha_2 > N$ and $\alpha_1 > N/2$.
Suppose that $\alpha_1 > \dots > \alpha_p > N/2 > \alpha_{p+1}$ and 
let $\beta$ be the sequence defined by
$$
\beta
 =
\begin{cases}
(N-\alpha_1, \dots, N-\alpha_p, \alpha_{p+1}, \dots, \alpha_r)
 &\text{if $p$ is even,} \\
(N-\alpha_1, \dots, N-\alpha_{p+1}, \alpha_{p+2}, \dots, \alpha_r)
 &\text{if $p$ is odd and $\alpha_p + \alpha_{p+1} \ge N+1$,} \\
(N-\alpha_1, \dots, N-\alpha_{p-1}, \alpha_p, \dots, \alpha_r)
 &\text{if $p$ is odd and $\alpha_p + \alpha_{p+1} \le N-1$.}
\end{cases}
$$
Since $\alpha_1 = \mu'_1 \le N$, we see that $\beta$ has at least two nonnegative entries.
However, if $\pi_{\Orth_N}( s_{[\mu]} ) = \pm S_{[(1^m)]}$, 
the rearrangement of $\beta$ in decreasing order is equal to 
$(m,-1, -2, \allowbreak \dots, \allowbreak -r+1)$, which has only one nonnegative entry.
This is a contradiction, so we have $\mu'_1 + \mu'_2 \le N$.
\end{demo}

By considering the tensor power of the defining representation and using Corollary~\ref{cor:Pieri-standard}, 
we can prove the following corollary, 
which is the special case $\alpha = (1, \dots, 1)$ of Theorem~\ref{thm:main2}. 

\begin{corollary}
\label{cor:Burrill}
\begin{enumerate}
\item[(1)]
For nonnegative integers $k$, $n$ and $m$, the following two sets are equinumerous:
\begin{enumerate}
\item[(a)]
The set of oscillating tableaux $(\lambda^{(i)})_{i=0}^k$ of length $k$ and shape $(m)$ 
such that $l(\lambda^{(i)}) \le n$ for each $i$.
\item[(b)]
The set of standard tableaux whose shape $\lambda$ satisfies $|\lambda| = k$, 
$l(\lambda) \le 2n$ and $c(\lambda) = m$.
\end{enumerate}
\item[(2)]
For nonnegative integers $k$, $N$ and $m$, the following two sets are equinumerous:
\begin{enumerate}
\item[(a)]
The set of oscillating tableaux $(\lambda^{(i)})_{i=0}^k$ of length $k$ and shape $(1^m)$ 
such that $\lambda^{(i)} \in \Par(\Orth_N)$ for each $i$.
\item[(b)]
The set of standard tableaux whose shape $\lambda$ satisfies $|\lambda| = k$, 
$l(\lambda) \le N$ and $r(\lambda) = m$.
\end{enumerate}
\item[(3)]
For nonnegative integers $k$, $N$ and $m$, the following two sets are equinumerous:
\begin{enumerate}
\item[(a)]
the set of sequences $(\lambda^{(i)})_{i=0}^k$ of partitions satisfying the following conditions:
\begin{enumerate}
\item[(i)]
$\lambda^{(0)} = \emptyset$, $\lambda^{(k)} = (1^m)$.
\item[(ii)]
$l(\lambda^{(i)}) \le n$ for each $i$.
\item[(iii)]
One of the following holds:
\begin{enumerate}
\item[(iii-1)]
$\lambda^{(i-1)} \subset \lambda^{(i)}$ and $|\lambda^{(i)}| = |\lambda^{(i-1)}| + 1$.
\item[(iii-1)]
$\lambda^{(i-1)} \supset \lambda^{(i)}$ and $|\lambda^{(i)}| = |\lambda^{(i-1)}| - 1$.
\item[(iii-1)]
$\lambda^{(i-1)} = \lambda^{(i)}$ and $l(\lambda^{(i-1)}) = n$.
\end{enumerate}
\end{enumerate}
\item[(b)]
The set of standard tableaux whose shape $\lambda$ satisfies $|\lambda| = k$, 
$l(\lambda) \le 2n+1$ and $r(\lambda) = m$ or $2n+1-m$.
\end{enumerate}
\item[(4)]
For nonnegative integers $k$, $n$ and $m$, the following two number coincide:
\begin{enumerate}
\item[(a)]
the summation
$$
\sum_O 2^{d(O)}
$$
over all oscillating tableaux $O = (\lambda^{(i)})_{i=0}^k$ of shape $(1^m)$ 
such that $l(\lambda^{(i)}) \le n$ for each $i$.
Here the statistic $d(O)$ is given by
$$
d(O) = \# \{ i : l(\lambda^{(i-1)}) = n, \ l(\lambda^{(i)}) < n \}.
$$
\item[(b)]
the number of standard tableaux whose shape $\lambda$ satisfies $l(\lambda) \le 2n$ 
and $r(\lambda) = m$ or $2n-m$.
\end{enumerate}
\end{enumerate}
\end{corollary}

It would be interesting to find bijective proofs of Theorem~\ref{thm:main2} 
and Corollary~\ref{cor:Burrill} by generalizing the arguments in \cite{BCFMM} and \cite{K}.


\end{document}